%% file: cauchy.tex
\documentclass[footinclude=false,11pt]{article}

\usepackage{fullpage}
\usepackage{amsmath}
\usepackage{amssymb}
\usepackage{xcolor}
\usepackage{graphicx}
\usepackage[export]{adjustbox}
\usepackage{tikz}
\usepackage{pgfplots}
\usepgfplotslibrary{patchplots}
\usetikzlibrary{shadows,backgrounds,arrows,positioning,matrix,calc}
\pgfplotsset{compat=newest}
\RequirePackage[tt=false]{libertine}
\RequirePackage[varqu]{zi4}
\RequirePackage[libertine]{newtxmath}
\RequirePackage[T1]{fontenc}
\usepackage[framemethod=TikZ]{mdframed}

\definecolor{myred}{rgb}{0.86,0.00,0.00}
\definecolor{myredlight}{rgb}{0.97,0.75,0.75}
\definecolor{myredlighter}{rgb}{0.99,0.94,0.94}
\definecolor{myredlighterr}{rgb}{1.0,0.98,0.98}
\definecolor{myblue}{rgb}{0.00,0.20,0.70}
\definecolor{mybluelight}{rgb}{0.75,0.80,0.93}
\definecolor{mybluelighter}{rgb}{0.94,0.95,0.98}
\definecolor{mybluelighterr}{rgb}{0.98,0.99,1.0}
\definecolor{mygreen}{rgb}{0.10,0.50,0.10}
\definecolor{mygreenlight}{rgb}{0.78,0.88,0.78}
\definecolor{mygreenlighter}{rgb}{0.94,0.97,0.94}
\definecolor{mygreenlighterr}{rgb}{0.99,0.99,0.99}
\definecolor{mygrey}{rgb}{0.40,0.40,0.40}
\definecolor{mygreylight}{rgb}{0.85,0.85,0.85}
\definecolor{mygreylighter}{rgb}{0.96,0.96,0.96}
\definecolor{mygreylighterr}{rgb}{0.99,0.99,0.99}
\definecolor{myorange}{rgb}{1.0,0.50,0.00}
\definecolor{myorangelight}{rgb}{1.0,0.87,0.75}
\definecolor{myorangelighter}{rgb}{1.0,0.96,0.93}
\definecolor{myorangelighterr}{rgb}{1.0,0.99,0.98}

\usepackage{listings}
\definecolor{lstshade}{gray}{0.95}
\definecolor{lstframe}{gray}{0.80}
\definecolor{lstcomment}{gray}{0.5}
\definecolor{lstattrib}{rgb}{0,0.34,0}
\lstnewenvironment{cpp}[1][]{\lstset{language=c++, #1}}
	{}
\lstnewenvironment{pseudocode}[1][]{\lstset{#1}}
	{}

\input{notations.tex}

\begin{document}

\title{
	Integration and Simulation of Bivariate Projective-Cauchy Distributions 
	within Arbitrary Polygonal Domains
}
\author{
	Jonathan Dupuy \hspace{1cm} Laurent Belcour \hspace{1cm} Eric Heitz
}
\date{\vspace{-0.8cm}}

\maketitle

\begin{abstract}
%It is known that the univariate Cauchy distribution arises from a projection 
%The univariate Cauchy distribution is defined on th
%It is known that the univariate Cauchy distribution arises from the stereographic 
%projection of a uniform variate on the unit circle. What is perhaps less known is  
%that this 

Consider a uniform variate on the unit upper-half sphere of dimension $d$. 
It is known that the straight-line projection through the center of the unit sphere 
onto the plane above it distributes this variate according to a $d$-dimensional
projective-Cauchy distribution.
\iffalse
Under this intuition, we show that integrating and simulating a univariate Cauchy 
distribution within an arbitrary line segment \mbox{$L \subset \mathbb R$} 
translates into respectively measuring and uniformly-sampling the angle subtended by 
this segment as seen from the origin of $\mathcal{S}^1$. We then extend this result to 2D. 
Specifically, we show that the same projection, which maps the hemisphere 
\mbox{$\mathcal{H}^2\subset \mathcal{S}^2$} 
onto the plane $\mathbb{R}^2$, links the solid-angle metric to the bivariate 
projective-Cauchy measure. 
\fi
In this work, we leverage the geometry of this construction in dimension $d=2$ to 
derive new properties for the bivariate projective-Cauchy distribution. 
Specifically, we reveal via geometric intuitions that integrating and 
simulating a bivariate projective-Cauchy distribution within an arbitrary domain 
translates into respectively measuring and sampling the solid angle subtended by the 
geometry of this domain as seen from the origin of the unit sphere.
To make this result practical for, e.g., generating truncated variants of the bivariate 
projective-Cauchy distribution, we extend it in two respects.
First, we provide a generalization to Cauchy distributions parameterized by 
location-scale-correlation coefficients. 
%We extend this result to Cauchy distributions parameterized by location-scale-correlation
%coefficients. 
Second, we provide a specialization to polygonal-domains, which leads to closed-form 
expressions. We provide a complete MATLAB implementation for the case of triangular 
domains, and briefly discuss the case of elliptical domains and how to further extend our 
results to bivariate Student distributions.
%In addition, we also discuss the case of elliptic domains.
%Our results can be used for, e.g., generating truncated variants of the bivariate 
%projective-Cauchy distribution. 
%A complete matlab implementation is provided 
%in appendix for triangular domains.

\end{abstract}

%\clearpage

% ----------------------------------------------------------------------------------------
% ----------------------------------------------------------------------------------------
\begin{figure}[t]
	\begin{center}
		\begin{tabular}{cc}
			\input{./fig_ctor_gnomonic_1d.tex}
			&
			\input{./fig_ctor_gnomonic_2d.tex}
			\vspace{-0.2cm}
			\\
			(1) & (2)
		\end{tabular}
	\end{center}
	\vspace{-0.6cm}
	\caption{Projective construction of the (1) univariate and 
	(2) bivariate Cauchy distribution.}
	\label{fig_ctor}
\end{figure}
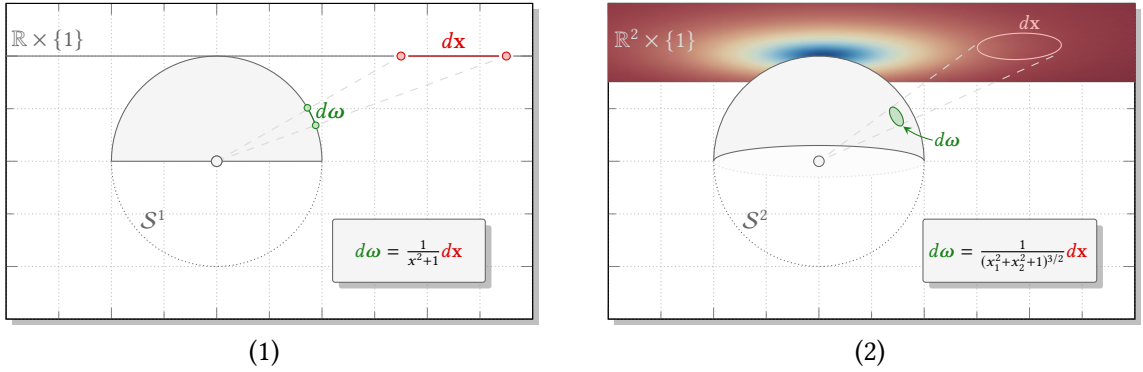
% ----------------------------------------------------------------------------------------
% ----------------------------------------------------------------------------------------

% ----------------------------------------------------------------------------------------
% ----------------------------------------------------------------------------------------
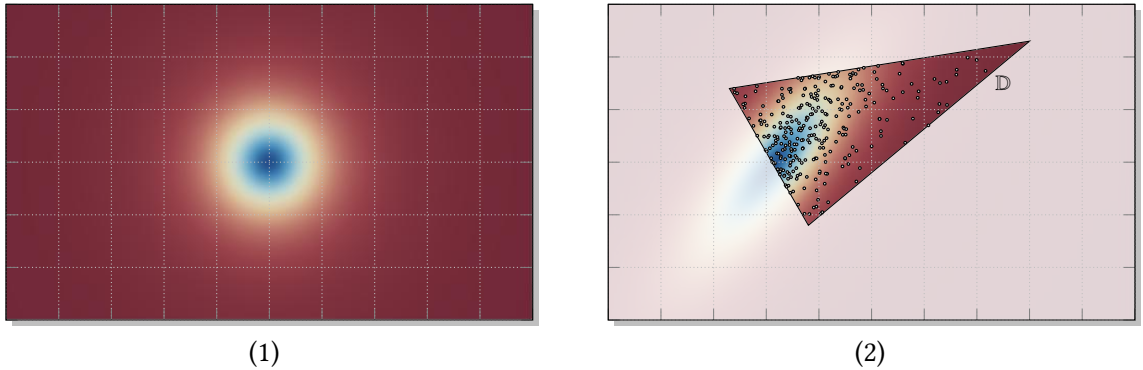
\begin{figure}[h]
	\begin{center}
		\begin{tabular}{cc}
			\input{./fig_pdf_std.tex}
			&
			\input{./fig_pdf_int.tex}
			\vspace{-0.2cm}
			\\
			(1) & (2)
		\end{tabular}
	\end{center}
	\vspace{-0.6cm}
	\caption{Isocontours of (1) the bivariate Cauchy density. We show how to 
	(2) integrate and simulate random variates distributed according to such a density 
	over an arbitrary domain $\mathbb{D}$. }
	\label{fig_intro}
\end{figure}
% ----------------------------------------------------------------------------------------
% ----------------------------------------------------------------------------------------

% ========================================================================================
% ========================================================================================
% ========================================================================================
\section{Introduction}

% ========================================================================================

% ========================================================================================
\paragraph{Construction of the Projective-Cauchy Distribution.}
A classic way of constructing a projective-Cauchy distribution is the 
following~\cite{knight1976caracterisation,dunau1988}:
Consider the projection of the upper-half unit 
sphere $\mathcal{S}^d$ of $\mathbb{R}^{d+1}$, $d \geq 1$, onto the plane 
$\mathbb{R}^{d} \times \{1\} \subset \mathbb{R}^{d+1}$ defined as
\begin{align}
	\label{eq_gnomonic}
	g_\perp(\bomega) 
	&= 
	(\omega_1/\omega_{d+1}, \,\ldots, \omega_{d}/\omega_{d+1}, \, 1), 
	%\\
	%\notag
	%\textrm{where } \bomega &= (\omega_1, \, \ldots,\, \omega_{n+1}) \in \mathcal{S}^n
	%\quad \textrm{having} \;\; \|(\omega_1, \, \ldots,\, \omega_{n+1}) \| = 1.
%	\\
%	\Rightarrow
%	g^{-1}(x_1, \, \ldots,\, x_{n}, 1) 
%	&= 
%	(x_1, \, \ldots,\, x_{n}, 1) / \sqrt{x_1^2 + \ldots + x_{n}^2 + 1}
\end{align}
where we use the notation $\bomega = (\omega_1, \, \ldots,\, \omega_{d+1}) \in \mathcal{S}^d$,
 i.e., $\|\bomega\| = 1$.
The application of such a projection to a uniform variate on the upper-half unit 
sphere $\mathcal{S}^d$ of $\mathbb{R}^{d+1}$ distributes it on the 
$\mathbb{R}^{d} \times \{1\}$ 
plane according to a standard projective-Cauchy distribution\footnote{
	Existing literature~\cite{knight1976caracterisation,dunau1988} consider a uniform 
	variate on the unit sphere $\mathcal S ^d$ rather than on its upper-half as we do here. 
	Since the projection defined in Equation~(\ref{eq_gnomonic}) is a two-to-one mapping 
	as for any \mbox{$\bomega \in \mathcal S ^d$}, we have \mbox{$g_\perp(\bomega) = g_\perp(-\bomega)$},
	it follows trivially that applying $g_\perp$ to a uniform variate on the upper-half unit 
	sphere of $\mathbb{R}^{d+1}$ also distributes it on the $\mathbb{R}^{d} \times \{1\}$ 
	plane according to a standard projective-Cauchy distribution.} 
with density
\begin{equation}
\label{eq_pdf_nd}
p_\std(\bx)
=
\frac{\Gamma(\alpha)}{\pi^{\alpha} }
\frac{1}{(x_1^2 + \cdots + x_d^2 + 1)^{\alpha}}, \quad \alpha=\frac{d+1}{2},
\end{equation}
where we use the notation $\bx = (x_1, \, \ldots, \, x_d, \, 1) \in \mathbb{R}^{d} \times \{1\}$;
Figure~\ref{fig_ctor} illustrates the geometry of this construction for the $d=\{1,2\}$ cases.

\iffalse
% ========================================================================================
\paragraph{Bijective-Projective Construction.}
The projection defined in Equation~(\ref{eq_gnomonic}) is a two-to-one mapping since for 
any \mbox{$\bomega \in \mathcal S ^d$}, we have \mbox{$g_\perp(\bomega) = g_\perp(-\bomega)$}.
It follows trivially that applying $g_\perp$ to a uniform variate on the upper-half unit 
sphere of $\mathbb{R}^{d+1}$ also distributes it on the $\mathbb{R}^{d} \times \{1\}$ 
plane according to a standard projective-Cauchy distribution; Figure~\ref{fig_ctor} 
illustrates the geometry of this construction for the $d=\{1,2\}$ cases.
In this work, we focus on this bijective-projective construction, which pictures 
$g_\perp$ as a bijection with inverse
\begin{equation}
	\label{eq_gnomonic_inv}
	g_\perp^{-1}(\bx) = 
	\frac{
		(x_1, \, \ldots, \, x_d, 1)
	}{\|(x_1, \, \ldots, \, x_d, 1)\|},
\end{equation}
%\iffalse
The bijective-projective also provides a geometric interpretation for the 
projective-Cauchy density, since the Equation~\ref{eq_pdf_nd} now translates into 
the Jacobian of the projection, i.e.,
\begin{equation}
	p_\std = \left \| \frac{dg_\perp^{-1}}{d\bx} \right \|.
\end{equation}
%\fi
\fi

% ========================================================================================
\paragraph{Bivariate Case.}
In this work, we focus on the specific case of the standard bivariate projective-Cauchy 
probability distribution, whose density is given by Equation~(\ref{eq_pdf_nd}) with $d=2$, 
i.e.,
\begin{equation}
	\label{eq_pdf_std}
	f_\std(\bx) = \frac{1}{2\pi} \frac{1}{(x_1^2 + x_2^2 + 1)^{3/2}}
	.
\end{equation}
%where the factor $1/2\pi$ corresponds to the probability density of the uniform 
%distribution on the hemisphere; 
Figure~\ref{fig_intro}~(1) shows a color-mapped plot of 
such a density, which is radially-symmetric and heavy-tailed. 
Hereafter, we refer to the standard bivariate projective-Cauchy distribution simply as 
the bivariate Cauchy distribution. 

% ========================================================================================
\paragraph{Contributions and Outline.}
In the following sections, we leverage the geometry of the projective construction 
to derive new properties for the bivariate Cauchy distribution.
Specifically, we show that the projective construction links 
the bivariate Cauchy distribution to the solid-angle metric in Section~\ref{sec_proof}.
Based on this result, we introduce a systematic approach 
to simulate and integrate the bivariate Cauchy distribution against an 
arbitrary domain $\mathbb D \subseteq \mathbb R ^2$ in 
Section~\ref{sec_integration_and_sampling}; this is useful, for, e.g., generating 
truncated variants of the bivariate Cauchy distribution.
Next, we extend our results to more general bivariate Cauchy densities parameterized 
by location-scale-correlation coefficients in Section~\ref{sec_generalization}. 
Finally, we consider the special case where the domain $\mathbb D$ represents 
the area enclosed by an arbitrary triangle in Section~\ref{sec_polygon}; 
Figure~\ref{fig_intro}~(2) shows an example of simulation of a location-scale-correlation 
Cauchy variate within a triangle using our contributions.  
A complete MATLAB implementation is provided in the supplemental material of this work.

%\clearpage

% ========================================================================================
% ========================================================================================
% ========================================================================================
\section{The Bivariate Cauchy Metric as the projection of Solid Angles}
\label{sec_proof}

In this section, we introduce and prove the main result of this paper. For this, we 
need to consider the mapping $g_\perp$ defined in Equation~(\ref{eq_gnomonic}) 
along with its inverse $g_\perp ^{-1}$, which is defined as
\begin{equation}
	\label{eq_gnomonic_inv}
	g_\perp^{-1}(\bx) = 
	\frac{
		(x_1, \, \ldots, \, x_d, 1)
	}{\|(x_1, \, \ldots, \, x_d, 1)\|}.
\end{equation}

% ========================================================================================
% ========================================================================================
\begin{center}
\fbox{%
\parbox{0.97\textwidth}{%
\vspace{-0.5cm}
\paragraph{Main Result.}

In the case $d=2$, the bivariate Cauchy density evaluates 
to the product of the uniform density on the upper-half sphere $\mathcal S^2$, 
times the Jacobian due to the change of measure from 
\mbox{$d\bomega \subset \mathcal{S}^2$} to \mbox{$d\bx \subset \mathbb{R}^{2} \times \{1\}$}:
\begin{equation}
	\label{eq_main_result}
	f_\std(\bx) = 
	\frac{1}{2\pi} 
	\left \| \frac{dg_\perp^{-1}(\bx)}{d\bx} \right \|.
\end{equation}
Put in other words, the bivariate Cauchy distribution is a solid-angle measure
expressed on the plane subtented by the action of $g_\perp$ on the upper-half unit sphere 
$\mathcal S^2$.
  }
}
\end{center}

%TODO: write the inverse in the main result section.

%The bivariate Cauchy density translates directly into the solid angle under the 
%action .

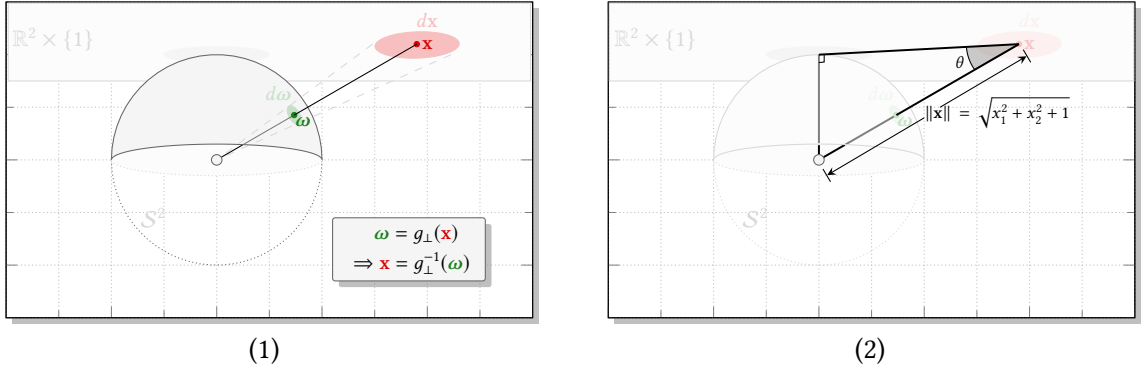
\begin{figure}[h]
	\begin{center}
		\begin{tabular}{cc}
		\input{./fig_ctor_gnomonic_2d_infinitesimals.tex}
		&
		\input{./fig_ctor_gnomonic_2d_triangle.tex}
		\vspace{-0.2cm}
		\\
		(1) & (2)
		\end{tabular}
	\end{center}
	\vspace{-0.6cm}
	\caption{Projective construction of the bivariate projective-Cauchy distribution.
	The projection operator maps (1) solid angles to areas, whose respective metrics
	are (2) linked through simple algebraic relations.}
	\label{fig_ctor_2d}
	\end{figure}

% ========================================================================================
% ========================================================================================
\paragraph{Proof.}
%We prove that the projective construction of the bivariate Cauchy distribution links 
%it to the solid angle metric. 
Let us consider an infinitesimal solid angle $d\bomega \subset \mathcal{S}^2$ 
centered around a directon $\bomega \in \mathcal S^2$ under the action of the projective 
mapping $g_\perp$. The projection maps $d\bomega$ to an infinitesimal area 
$d \bx$ located on the $\mathbb R^2 \cup \{1\}$ plane around the point 
$\bx = g_\perp(\bomega) \in \mathbb R^2 \cup \{1\}$; Figure~\ref{fig_ctor_2d}~(1) 
illustrates this setup.

% ----------------------------------------------------------------------------------------
% ----------------------------------------------------------------------------------------
% ----------------------------------------------------------------------------------------
% ----------------------------------------------------------------------------------------

By construction, the infinitesimal area $d \bx$ is located at distance 
$\|\bx\|= \|g_\perp(\bomega)\|$ from the origin and projected at an angle $\theta$ onto
the $\mathbb R^2 \cup \{1\}$ plane; Figure~\ref{fig_ctor_2d}~(2) illustrates this 
geometric relationship. Since solid angles scale quadratically with distance and 
oriented projections produces a foreshortening effect of $\sin \theta$, we have
\begin{align}
	\label{eq_forshortening}
	\|\bx\|^2 \, d\bomega 
	&= \sin \theta \, d \bx.
\end{align}
Now, using the identities $\|\bx\| = \sqrt{x_1^2 + x_2^2 + 1}$ and $\sin \theta = 1/\|\bx\|$, 
we further obtain 
\begin{equation}
	\label{eq_jacobian_2d}
  d\bomega = \frac{1}{(x_1^2 + x_2^2 + 1)^{3/2}} \, d \bx.
\end{equation}
Hence, the Jacobian $\left\|\frac{d\bomega}{d \bx}\right \|$ is directly proportional to 
Equation~(\ref{eq_pdf_std}). 
Under the projective construction, the algebraic expression of the bivariate 
Cauchy density translates into the Jacobian of the projection, which is 
Equation~(\ref{eq_jacobian_2d}), normalized by the solid angle of the hemisphere, 
which is $2 \pi$. Using the equality $\bomega = g_\perp^{-1}(\bx)$, we finally obtain 
our main result: 
\begin{equation}
f_\std( \bx )
= \frac{1}{2\pi}\left\|\frac{d \bomega}{d \bx} \right\|\\
= \frac{1}{2\pi}\left\|\frac{d g_\perp^{-1}(\bx)}{d \bx} \right\|.
\end{equation}
This concludes our proof.

%\clearpage
% ========================================================================================
% ========================================================================================
% ========================================================================================
\section{Integration and Simulation}
\label{sec_integration_and_sampling}

In this section, we introduce a systematic approach for simulating and integrating 
the bivariate Cauchy distribution within an arbitrary domain 
$\mathbb D \subseteq \mathbb R ^2 \times \{1\}$. This requires the following corollary result: 

% ========================================================================================
\begin{center}
	\fbox{%
	\parbox{0.97\textwidth}{%
	\vspace{-0.5cm}
\paragraph{Corollary Result.}
%Consider an arbitrary domain $\mathbb D \subset \mathbb R ^2$. A restricted measure 
%restricted to $\mathbb D$  
Consider an arbitrary domain $\mathbb D \subseteq \mathbb R ^2 \times \{1\}$. 
The inverse projection $g_\perp^{-1}$ maps $\mathbb D$ to the solid angle it subtends 
with respect to the unit sphere $\mathcal{S}^2$ of $\mathbb{R}^{3}$.
	}
}
\end{center}

\begin{figure}[h]
	\begin{center}
		\begin{tabular}{cc}
		\input{./fig_simulate.tex}
		&
		\input{./fig_integrate.tex}
		\vspace{-0.2cm}
		\\
		(1) & (2)
		\end{tabular}
	\end{center}
	\vspace{-0.6cm}
	\caption{Geometric interpretation of our systematic methodology for (1) simulating 
	and (2) integrating bivariate Cauchy distributions within an arbitrary domain 
	$\mathbb D \subseteq \mathbb R ^2 \times \{1\}$. We re-express both operations on 
	the surface of the unit sphere $\mathcal S ^2$.}
	\label{fig_sim_int}
	\end{figure}
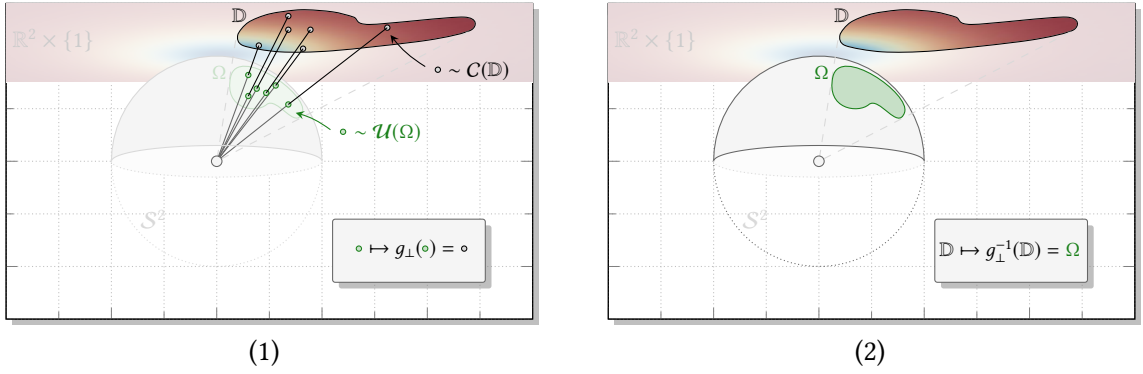

% ========================================================================================
\paragraph{Variate Simulation.}
The projective construction of the bivariate Cauchy distribution implies a 
straightforward simulation mechanism: it suffices to generate a uniform variate 
within the solid angle \mbox{$\Omega = g_\perp^{-1}(\mathbb D)$} and project it using 
$g_\perp$, i.e.,
Equation~(\ref{eq_gnomonic}) with $d=2$; Figure~\ref{fig_sim_int}~(1) illustrates this 
approach.
It follows that if the geometry of $\Omega$ can be analytically sampled with uniform 
probability, then the bivariate Cauchy distribution can be simulated analytically 
within the area enclosed by $\mathbb D$.

% ========================================================================================
\paragraph{Bounded Integration.}
The projective construction of the bivariate Cauchy distribution implies a 
straightforward method for solving integrals of the form
\begin{align}
	\label{eq_integral_std}
	\mathcal{I} = 
	\int_\mathbb{D}
	f_\std( \bx ) \, d \bx.
\end{align}
Let $\Omega$ denote the solid angle subtented by $\mathbb D$ with respect to the unit 
sphere unit sphere $\mathcal{S}^2$ of $\mathbb{R}^{3}$, i.e., 
\mbox{$\Omega = g_\perp^{-1}(\mathbb D)$}.
We re-express Equation~(\ref{eq_integral_std}) as a measure over solid-angles using 
the substitution \mbox{$\bx = g_\perp(\bomega)$}, which makes integrand constant
\begin{equation}
	\label{eq_integral_std_result}
	\mathcal{I} 
	=
	\int_{\Omega} 
	\underbrace{f_\std( g_\perp(\bomega))}_{=\frac{1}{2\pi}
	\left\|\frac{d \bomega}{d \bx}\right\|, \, \textrm{Eq.(\ref{eq_main_result}})} 
	\left\|\frac{d \bx}{d \bomega}\right\| 
	d\bomega
	= \int_{\Omega} \frac{1}{2\pi} \; d\bomega, \quad \Omega = g_\perp^{-1}(\mathbb D).
\end{equation}
Put in other words, integrating the bivariate Cauchy distribution reduces to computing 
a solid angle; Figure~\ref{fig_sim_int}~(2) illustrates this approach.
It follows that if the solid angle of $\Omega$ can be computed analytically, then 
the bivariate Cauchy distribution can be integrated analytically within the area 
enclosed by $\mathbb D$. 

% ========================================================================================
%\paragraph{Truncated Distribution.}

\clearpage
% ----------------------------------------------------------------------------------------
% ----------------------------------------------------------------------------------------
% ----------------------------------------------------------------------------------------
% ----------------------------------------------------------------------------------------

\vspace{-0.5cm}
% ========================================================================================
% ========================================================================================
% ========================================================================================
\section{Generalization to Location-Scale-Correlation}
\label{sec_generalization}

In this section, we generalize the results of the previous section to a more 
general parameterization of the bivariate Cauchy distribution.

\begin{figure}[h]
	\begin{center}
		\begin{tabular}{cc}
		\input{./fig_pdf_std.tex}
		&
		\input{./fig_pdf.tex}
		\vspace{-0.2cm}
		\\
		(1) & (2)
		\end{tabular}
	\end{center}
	\vspace{-0.6cm}
	\caption{Isocontours of (1) the standard projective-Cauchy density, and 
	(2) a general projective-Cauchy density with location parameters $a_1 = -1.9, a_2 = -0.1$, 
	scale parameters $b_1 = 1.4, b_2 = 1.7$, and correlation $\rho = 0.8$.}
	\label{fig_pdf}
\end{figure}
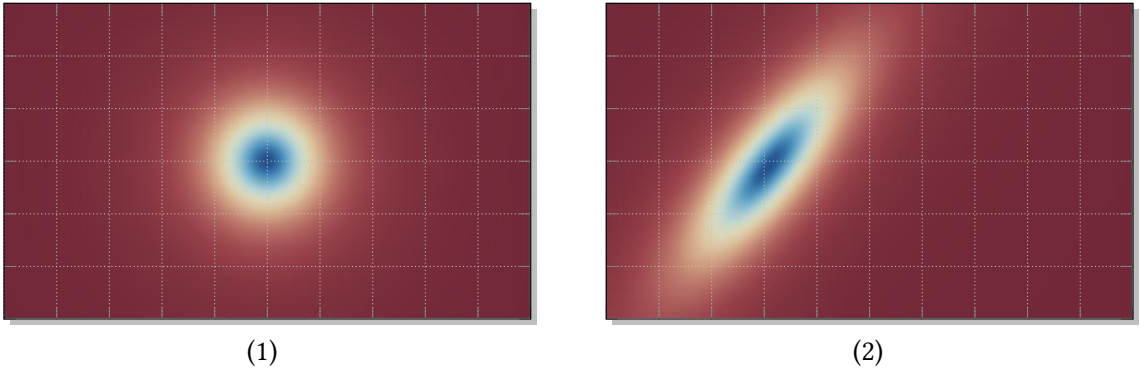

% ========================================================================================
\paragraph{Parameterization.}
The bivariate Cauchy density exhibits radially symmetric isocontours.
These can be extended to more general elliptic isocontours when the random 
variates are transformed according to a location-scale-correlation substitution. 
%parameterize radially symmetric random variates is to warp them via 
%the location-scale-correlation transformation. 
For the 2D case, considering the location parameters \mbox{$a_1, a_2 \in \mathbb{R}$}, 
scale parameters \mbox{$b_1, b_2 > 0$}, and correlation parameter \mbox{$\rho \in (-1,1)$}, 
the transformation takes the form 
\begin{equation}
	\label{eq_warp}
	g_{\mathcal E}(\bx) = \left(
		b_1 \, x_1 + a_1, \;
		b_2 \left[ \rho x_1 + x_2 {\textstyle\sqrt{1-\rho^2}} \right] + a_2,
		1
	\right).
\end{equation}
The application of $g_{\mathcal E}$ on a random variate distributed as a bivariate 
Cauchy distribution redistributes it according to a more general bivariate Cauchy 
distribution with density
\begin{align}
	\label{eq_pdf}
	f(\bx ;  a_1, a_2, b_1, b_2, \rho)
	%f(x_1, x_2)
	&= f_\std
	\left(
		g^{-1}_{\mathcal E}(\bx)
	\right) 
	\left\|\frac{dg^{-1}_{\mathcal E}(\bx)}{d \bx}\right\|, \\
	\label{eq_jacobian}
	\left\|\frac{dg^{-1}_{\mathcal E}(\bx)}{d \bx}\right\| &= \frac{1}{b_1 b_2 \sqrt{1 - \rho^2}},
\end{align}
where $g^{-1}_{\mathcal E}$ denotes the inverse of $g_{\mathcal E}$, which takes the form
\begin{equation}
	\label{eq_warp_inverse}
	g^{-1}_{\mathcal E}(\bx) = \left(
		\frac{x_1 - a_1}{b_1},
		\frac{b_1(x_2 - a_2) - \rho \, b_2 (x_1 - a_1)}{b_1 b_2 \sqrt{1-\rho^2}},
		1
	\right) .
\end{equation}
The standard bivariate Cauchy case corresponds to setting the parameters to 
$a_1=a_2=0$, $b_1=b_2=1$, and $\rho=0$. 
Figure~\ref{fig_pdf} plots the isocontours of the standard bivariate Cauchy density 
against that of an elliptic one. Hereafter, we refer to such generalized bivariate 
Cauchy distributions as elliptical bivariate Cauchy distributions.
Note that Equation~(\ref{eq_pdf}) is equivalent to the more standard expression
\begin{equation}
	f(\bx ;  a_1, a_2, b_1, b_2, \rho)
	=
	\frac{1}{2 \pi b_1 b_2 \sqrt{1 - \rho^2}} 
	\left( 1 + \frac{z}{1-\rho^2} \right)^{-\frac{3}{2}},
\end{equation}
where
\begin{equation}
z = 
\frac{(x_1-a_1)^2}{b_1^2} 
+ \frac{(x_2-a_2)^2}{b_2^2} 
- \frac{2\rho(x_1-a_1)(x_2-a_2)}{b_1 b_2}.
\end{equation}

% ----------------------------------------------------------------------------------------
% ----------------------------------------------------------------------------------------
% ----------------------------------------------------------------------------------------
% ----------------------------------------------------------------------------------------
\vspace{-0.5cm}

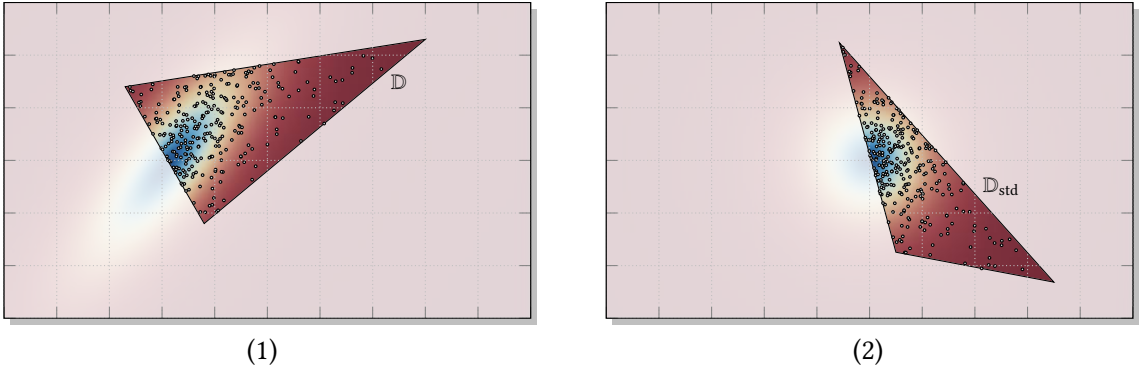
\begin{figure}[h]
	\begin{center}
		\begin{tabular}{cc}
		\input{./fig_pdf_int.tex}
		&
		\input{./fig_pdf_std_int.tex}
		\vspace{-0.2cm}
		\\
		(1) & (2)
		\end{tabular}
	\end{center}
	\vspace{-0.6cm}
	\caption{Invariance of polygonal integrals over Cauchy densities. (1) The integral 
	of any elliptical bivariate Cauchy density against a domain $\mathbb D$ is equal to 
	(2) another integral of a bivariate Cauchy density over a domain $\mathbb D_\std$.}
	\label{fig_integral}
\end{figure}

% ========================================================================================
\paragraph{Variate Simulation.}
The location-scale-correlation construction of the elliptical bivariate Cauchy 
distribution implies a straightforward simulation mechanism within an arbitrary domain 
$\mathbb D \subset \mathbb R ^2 \times \{1\}$: it suffices to generate a bivariate Cauchy 
variate within the domain \mbox{$\mathbb D_\std = g_{\mathcal E}^{-1}(\mathbb D)$} using the 
results from the previous section, i.e., uniformily sample the spherical 
domain $\Omega = g_\perp^{-1}(\mathbb D_\std)$, and then transform it using 
Equation~(\ref{eq_warp}); Figure~\ref{fig_integral} illustrates this approach.
It follows that if the geometry of $\Omega = g_\perp^{-1}(\mathbb D_\std)$ can be 
analytically sampled with uniform probability, then the bivariate Cauchy distribution 
can be simulated analytically within the area enclosed by $\mathbb D$.  

% ========================================================================================
\paragraph{Bounded Integration.}
The location-scale-correlation construction of the elliptical bivariate Cauchy 
distribution implies a straightforward method for solving integrals of the form
\begin{align}
	\label{eq_integral_elliptical}
	\mathcal{J} = 
	\int_\mathbb{D}
	f(\by ;  a_1, a_2, b_1, b_2, \rho) \, d\by.
\end{align}
Let $\mathbb D_\std$ denote the domain resulting from the application of $g_{\mathcal E}$ 
to $\mathbb D$. We re-express Equation~(\ref{eq_integral_elliptical}) as a measure over 
the bivariate Cauchy distribution using the substitution $\by = g_{\mathcal E}(\bx)$, 
i.e.,
\begin{equation}
	\mathcal{J}
	= 
	\int_{\mathbb D_\std} 
	\underbrace{
		f(g_{\mathcal E}(\bx); a_1, a_2, b_1, b_2, \rho)
	}_{
		= f_\std(\bx) \, \left \|\frac{d \bx}{d \by}\right\| \textrm{, Eq.(\ref{eq_pdf})}
	}
	\,
	\left \|\frac{d \by}{d \bx}\right\| \, d\bx, 
	= 
	\int_{\mathbb D_\std} f_\std(\bx) \, d\bx, 
	\quad 
	\mathbb D_\std = g_{\mathcal E}^{-1}(\mathbb D).
\end{equation}
Using the results of the previous section, we re-express this new integral as a uniform 
measure over solid-angles, where the integrand becomes constant
\begin{equation}
	\mathcal{J}
	= 
	\int_{\Omega} 
	\underbrace{f_\std( g_\perp(\bomega))}_{=\frac{1}{2\pi}
	\left\|\frac{d \bomega}{d \bx}\right\|, \, \textrm{Eq.(\ref{eq_main_result}})} 
	\left\|\frac{d \bx}{d \bomega}\right\| 
	d\bomega
	=
	\int_{\Omega} \frac{1}{2\pi} \;d\bomega,
	\quad 
	\Omega = g_{\perp}^{-1}(\mathbb D_\std).
\end{equation}
It follows that if the solid angle of $\Omega$ can be computed analytically, then 
the elliptical bivariate Cauchy distribution can be integrated analytically within 
the area enclosed by $\mathbb D$.

% ----------------------------------------------------------------------------------------
% ----------------------------------------------------------------------------------------
% ----------------------------------------------------------------------------------------
% ----------------------------------------------------------------------------------------
\vspace{-0.5cm}

%\clearpage
% ========================================================================================
% ========================================================================================
% ========================================================================================
\section{Polygonal Bounds}
\label{sec_polygon}

In this section, we specialize our results from Section~\ref{sec_integration_and_sampling} 
to the case where $\mathbb D$ represents the area enclosed by an arbitrary polygon. 
We also show that these results extend trivially to the 
location-scale-correlation generalization introduced in Section~\ref{sec_generalization}.

\iffalse
Since the location-scale-correlation substitution is 
a linear transformation of the $\mathbb{R}^2$ plane, it acts as an automorphism on any 
polygons. Put in other words, if we can integrate and simulate the bivariate Cauchy density 
against arbitrary polygons, then we can also integrate arbitrary 
location-scale-correlation densities against arbitrary polygons
\fi

% ========================================================================================
\paragraph{Variate Simulation.}
Simulating a bivariate Cauchy distribution withing the area enclosed by a polygon 
requires uniformly sampling the solid angle of the spherical geometry it subtends via 
$g_\perp^{-1}$, which is a spherical polygon, 
and then projecting these samples using Equation~(\ref{eq_gnomonic}). 
If \mbox{$\bv_1, \ldots, \bv_n \in \mathbb R^2 \times \{1\}$} denote the vertices of the 
polygon, then the spherical polygon has vertices 
\mbox{$g_\perp^{-1}(\bv_1), \ldots, g_\perp^{-1}(\bv_n) \in \mathcal S^2$}.
For details on how to uniformly sampling such geometries, we refer the reader to the 
literature on the subject~\cite{arvo1995stratified,urena2013area}, as this is out of 
the scope of this paper; Appendix~\ref{sec_appendix} provides pseudocode for 
the case where the polygonal-domain is defined by an arbitrary triangle, and 
Figure~\ref{fig_integral}~(2) shows an example of variate simulation within a triangle. 

% ========================================================================================
\paragraph{Polygonal Integration.}
Integrating a bivariate Cauchy distribution against the area enclosed by a polygon 
requires computing the solid angle of the spherical geometry it subtends via 
$g_\perp^{-1}$, which is a spherical polygon. For a polygon with $N \geq 3$ 
vertices, the solid angle of the spherical polygon is given by
\begin{equation}
	\label{eq_solid_angle}
	\mathcal I = \sum_{n=1}^{N} A_n - (N - 2) \pi,
\end{equation}
where $A_n:= A_n(\Omega)$ denotes the $n$-th interior angle of the spherical polygon.
%If \mbox{$\bv_1, \ldots, \bv_n \in \mathbb R^2 \times \{1\}$} denote the vertices of the 
%polygon, then the spherical polygon has vertices 
%\mbox{$g_\perp^{-1}(\bv_1), \ldots, g_\perp^{-1}(\bv_n) \in \mathcal S^2$}.
%It follows trivially that the $n$-th interior angle of $\Omega$ is 
%\begin{equation}
%A_n = \cos^{-1}(\bomega_n \cdot \bomega_{n+1}), \quad \bomega_i = g_\perp^{-1}(\bv_i).
%\end{equation}
Appendix~\ref{sec_appendix} provides pseudocode for the case where the polygonal-domain 
is defined by an arbitrary triangle.

% ========================================================================================
\paragraph{Location-Scale-Correlation Generalization.}
The results derived here for the bivariate Cauchy distribution extend 
trivially to the location-scale-correlation generalization of 
Section~\ref{sec_generalization}: Since the mapping $g_{\mathcal E}$ defined in 
Equation~(\ref{eq_warp_inverse}) is a linear transformation, it acts as an automorphism 
on polygons. Put in other words, if we can integrate and simulate the 
bivariate Cauchy density against arbitrary polygons, then we can also integrate and 
simulate arbitrary location-scale-correlation densities against arbitrary polygons 
using the simulation and integration methods introduced in the previous section.
Appendix~\ref{sec_appendix} provides pseudocode for the case where the polygonal-domain 
is defined by an arbitrary triangle.

%\iffalse
%\clearpage
% ========================================================================================
% ========================================================================================
% ========================================================================================
\section{Discussion}
\label{sec_discussion}

In this section, we informally discuss the extension of our results to 
elliptical domains, and Student distributions with integer-valued degrees of freedom.
We defer rigorous proofs from our claims for future work.

% ========================================================================================
\paragraph{On Elliptical Domains.}
Following the methodology introduced in Section~\ref{sec_integration_and_sampling}, 
we can integrate and simulate the bivariate Cauchy density against 
arbitrary ellipses as long as we can deal with the solid angle they subtend 
on $\mathcal S ^2$.
Such geometries have been studied in the literature, and there exists methods to 
compute their solid angles~\cite{heitz2017computing}, as well as uniform-sample 
them~\cite{guillen2017area}.
Furthermore, the mapping $g_{\mathcal E}$ acts as an automorphism on ellipses, 
which means that we can also integrate and simulate arbitrary location-scale-correlation 
densities against arbitrary ellipses.

% ========================================================================================
\paragraph{On Student Distributions.}
As demonstrated by Equation~(\ref{eq_main_result}), the bivariate Cauchy distribution 
arises from the uniform distribution on the upper-half sphere $\mathcal S^2$ times 
the Jacobian of the projection $g_\perp$. This suggests a mechanism for deriving other 
distributions than just the Cauchy distribution by weighting directional 
distributions that differ from the uniform one. Interestingly, if we consider the 
power-law distribution $\omega_3^\nu$, $\nu \geq 1$, then using the identity due to 
Equation~(\ref{eq_gnomonic_inv})
\begin{equation}
\omega_3 = \frac{1}{\sqrt{x_1^2 + x_2^2 + 1}},
\end{equation}
we get  
\begin{equation}
d\bomega = \frac{1}{(x_1^2+x_2^2+1)^{\frac{2+\nu}{2}}} d\bx,
\end{equation}
which implies that the Jacobian $\left\|\frac{d\bomega}{d \bx}\right \|$ is now directly 
proportional to the bivariate Student distribution with $\nu$ degrees of freedom. 
As such, any integral of the form 
\begin{equation}
	\mathcal I' = \int_{\mathbb D} \frac{1}{(x_1^2+x_2^2+1)^{\frac{2+\nu}{2}}} d\bx,
\end{equation}
may be re-expressed in the directional domain as a 
\begin{equation}
	\mathcal I' = \int_{\Omega} \omega_3^\nu \, d\bomega, \quad \Omega = g_\perp^{-1}(\mathbb D).
\end{equation}
Since it is known that such integrals have analytical 
solutions over polygonal domains for any $\nu \in \mathbb N$~\cite{arvo1995applications}, 
we believe that new properties for the bivariate Student distribution may be derived 
based on our results for the bivariate Cauchy distribution.

\clearpage
% ========================================================================================
% ========================================================================================
% ========================================================================================
\section{Appendix: Pseudocode}
\label{sec_appendix}

% ----------------------------------------------------------------------------------------
% ----------------------------------------------------------------------------------------
\begin{pseudocode}[caption = Simulating and integrating a bivariate Cauchy density against a polygon.]
/* Project to the unit hemisphere */
function PlaneToHemisphere($\bx = \{x_1, \; x_2, \; 1\}$)
	return $\bomega = \bx \;/ \; \|\bx\|$

/* Project to the plane */
function HemisphereToPlane($\bomega = \{\omega_1, \; \omega_2, \; \omega_3\}$)
	return $\bx = \bomega \;/ \;  \omega_3$

/* Polygonal simulation of a standard bivariate Cauchy distribution */	
function SimulateCauchyStd($\mathcal{P} = \{\bv_1, \ldots, \bv_N\}$, $\bu = \{u_1, \; u_2\}$)
	$\Omega \; \leftarrow \;\{$PlaneToHemisphere($\bv_1$)$,\ldots , $ PlaneToHemisphere($\bv_N$)$\}$
	return HemisphereToPlane(SampleSolidAngle($\Omega$, $\bu$))

/* Polygonal integral of a standard bivariate Cauchy distribution */	
function IntegrateCauchyStd($\mathcal{P} = \{\bv_1, \ldots, \bv_N\}$)
	$\Omega \; \leftarrow \;\{$PlaneToHemisphere($\bv_1$)$,\ldots , $ PlaneToHemisphere($\bv_N$)$\}$
	return ComputeSolidAngle($\Omega$)
\end{pseudocode}
% ----------------------------------------------------------------------------------------
% ----------------------------------------------------------------------------------------
\vspace{0.5cm}
% ----------------------------------------------------------------------------------------
% ----------------------------------------------------------------------------------------
\begin{pseudocode}[caption = Simulating and integrating a bivariate Cauchy density against a polygon with location-scale-correlation parameters.]
/* LocationScaleCorrelationForward */
function LSCForward($\bx = \{x_1, \; x_2, \; 1\}$, $a_1,\; a_2,\; b_1,\; b_2,\; \rho$)
	/* See Equation (6) */
	$x_1' \; \leftarrow \; b_1 x_1 + a_1$
	$x_2' \; \leftarrow \; b_2 [ \rho x_1 + x_2 {\textstyle\sqrt{1-\rho^2}} ] + a_2$
	return $\{x_1', \; x_2', \; 1\}$

/* LocationScaleCorrelationForward */
function LSCBackward($\bx = \{x_1, \; x_2, \; 1\}$, $a_1,\; a_2,\; b_1,\; b_2,\; \rho$)
	/* See Equation (6) */
	$x_1' \; \leftarrow \; (x_1 - a_1) \;/\; b_1$
	$x_2' \; \leftarrow \; (b_1(x_2 - a_2) - \rho \, b_2 (x_1 - a_1))\;/\;(b_1 b_2 \sqrt{1-\rho^2})$
	return $\{x_1', \; x_2', \; 1\}$

/* Polygonal simulation of a standard bivariate Cauchy distribution */	
function SimulateCauchyElliptic($\mathcal{P} = \{\bv_1, \ldots, \bv_N\}$, $\bu = \{u_1, \; u_2\}$)
	return LSCForward(SimulateCauchy(LSCBackward($\mathcal{P}$), $\bu$))

/* Polygonal integral of a standard bivariate Cauchy distribution */	
function IntegrateCauchyElliptic($\mathcal{P} = \{\bv_1, \ldots, \bv_N\}$)
	return IntegrateCauchyStd(LSCBackward($\mathcal{P}$))
\end{pseudocode}
% ----------------------------------------------------------------------------------------
% ----------------------------------------------------------------------------------------

\iffalse
% ----------------------------------------------------------------------------------------
% ----------------------------------------------------------------------------------------
\begin{pseudocode}[caption = Integrating a bivariate Cauchy density against a polygon.]
/* Location-scale-correlation substitution */
function LocationScaleCorrelation($\bv = \{x_v, \; y_v\}, \,a_1, \,a_2, \,b_1, \,b_2, \,\rho$)
	/* See Equation (6) */
	$x_v' \; \leftarrow \; b_1 x_v + a_1$
	$y_v' \; \leftarrow \; b_2 [ \rho x_v + y_v {\textstyle\sqrt{1-\rho^2}} ] + a_2$
	return $\{x_v', \; y_v'\}$

/* Polygonal Integral of a general bivariate Cauchy distribution */	
function CauchyIntegral($\mathcal{P} = \{\bv_1, \ldots, \bv_N\}, \,a_1, \,a_2, \,b_1, \,b_2, \,\rho$)
	$\mathcal{P}' \; \leftarrow \;\{$Warp($\bv_1, \,a_1, \,a_2, \,b_1, \,b_2, \,\rho$)$, \ldots ,$ Warp($\bv_N, \,a_1, \,a_2, \,b_1, \,b_2, \,\rho$)$\}$
	return CauchyIntegralStd($\mathcal{P}'$)
\end{pseudocode}
% ----------------------------------------------------------------------------------------
% ----------------------------------------------------------------------------------------
\fi	
\clearpage

\bibliography{cauchy}{}
\bibliographystyle{plain}
\end{document}

%% file: notations.tex
\newcommand{\br}{\mathbf{r}}

\newcommand{\bu}{\mathbf{u}}
\newcommand{\bv}{\mathbf{v}}
\newcommand{\bx}{\mathbf{x}}
\newcommand{\by}{\mathbf{y}}

\newcommand{\bomega}{\boldsymbol{\omega}}

\newcommand{\std}{\textrm{std}}

%% file: fig_ctor_gnomonic_1d.tex
\pgfplotsset{every axis/.append style = {
	x = 1.4cm,
	y = 1.4cm,
	xmin = -2, xmax = 3,
	ymin = -1.5, ymax = 1.5,
	grid = both,
	axis line style = {color = black},
	minor grid style={densely dotted},
	major grid style={densely dotted},
	y label style = {rotate = -90},
	minor tick num = 1,
	xtick = {-3,...,3},
	ytick = {-2,...,2},
	xticklabels = {},
	yticklabels = {},
	no markers
}}
\begin{tikzpicture}
    \begin{axis}[name = myaxis]
		\addplot[samples=45, domain=pi:2*pi, color = mygrey, densely dotted] (
			{cos(deg(x))}, {sin(deg(x))}
		);
		\addplot[samples=45, domain=0:pi, color = mygrey, fill = mygreylighter] (
			{cos(deg(x))}, {sin(deg(x))}
		);
		\addplot[samples=45, domain=20:30, color = mygreen] (
			{cos(x)}, {sin(x)}
		);
		\draw (axis cs:-1,0) -- (axis cs:1,0) [mygrey];
		% real line
		\draw (axis cs:-3,1) -- (axis cs:3,1) [mygrey];
		% infinitesimal on the real line
		\draw (axis cs:1.82,1) -- (axis cs:2.68,1) [myred];
		% first point
		\draw (axis cs:0,0) -- (axis cs:1.75,1) [mygreylight, dashed];
		\draw[white, fill=white] (axis cs:1.75,1) circle [radius=0.125*0.5];
		\draw[myred, fill=myredlight] (axis cs:1.75,1) circle [radius=0.125*0.3];
		\draw[mygreen, fill=mygreenlight] (axis cs:0.861538,0.507692) circle [radius=0.125*0.25];
		% second point
		\draw (axis cs:0,0) -- (axis cs:2.75,1) [mygreylight, dashed];
		\draw[white, fill=white] (axis cs:2.75,1) circle [radius=0.125*0.5];
		\draw[myred, fill=myredlight] (axis cs:2.75,1) circle [radius=0.125*0.3];
		\draw[mygreen, fill=mygreenlight] (axis cs:0.939793, 0.341743) circle [radius=0.125*0.25];
		% projection center
		\draw[mygrey, fill=mygreylighter] (axis cs:0,0) circle [radius=0.125*0.4];
		% labels
		\node[myred, draw = none, anchor = center] at (axis cs:2.25, 1.15)
			{\scalebox{0.8}{$d\bx$}};
		\node[mygreen, draw = none, anchor = center] at (axis cs:1.08, 0.49)
			{\scalebox{0.8}{$d \bomega$}};
		% angles
		%\draw[mygrey, fill=mygreylighter] (axis cs:-1,0) circle [radius=0.125*0.2];
		%\draw[mygrey, fill=mygreylighter] (axis cs:+1,0) circle [radius=0.125*0.2];
		%\draw[mygrey, fill=mygreylighter] (axis cs:0,1) circle [radius=0.125*0.2];
		\draw[mygrey, fill=mygreylighter] (axis cs:0,0) circle [radius=0.125*0.4];
		%\node[mygrey, draw = none, anchor = center] at (axis cs:-0.8, -0.10)
		%	{\scalebox{0.5}{$\theta = \frac{\pi}{2}$}};
		%\node[mygrey, draw = none, anchor = center] at (axis cs:+1.23, -0.10)
		%	{\scalebox{0.5}{$\theta = -\frac{\pi}{2}$}};
		%\node[mygrey, draw = none, anchor = center] at (axis cs:+0.17, +1.12)
		%	{\scalebox{0.5}{$\theta = 0$}};
		% labels
		\node[mygrey, draw = none, anchor = center] at (axis cs:-0.6, -0.55)
			{\scalebox{0.8}{$\mathcal S ^1$}};
		\node[mygrey, draw = none, anchor = center] at (axis cs:-1.6, 1.15)
			{\scalebox{0.8}{$\mathbb R \times \{1\}$}};
		% formulas
		\draw[xshift = 4.0, yshift = -2.0, rounded corners=1, mygrey, fill = mygreylighter, drop shadow] 
		(1, -1.1) rectangle (2.45, -0.5);
		\node[xshift = 4.0, yshift = -2.0, black, draw = none, anchor = north west] at (1.07, -0.60)
			{\scalebox{0.7}{
				${\color{mygreen} d\bomega} =  \frac{1}{x^2+1} {\color{myred} d\bx}$
			}};
    \end{axis}
    \begin{pgfonlayer}{background}
    \draw[fill=white, drop shadow]
        (myaxis.north west) 
        rectangle (myaxis.south east);
    \end{pgfonlayer}
\end{tikzpicture}

%% file: fig_ctor_gnomonic_2d.tex
\pgfplotsset{every axis/.append style = {
	x = 1.4cm,
	y = 1.4cm,
	xmin = -2, xmax = 3,
	ymin = -1.5, ymax = 1.5,
	grid = both,
	axis line style = {color = black},
	minor grid style={densely dotted},
	major grid style={densely dotted},
	y label style = {rotate = -90},
	minor tick num = 1,
	xtick = {-3,...,3},
	ytick = {-2,...,2},
	xticklabels = {},
	yticklabels = {},
	no markers
}}
\begin{tikzpicture}
    \begin{axis}[name = myaxis]
		% real plane
		\draw[mygreylight, fill = mygreylighterr] 
			(-1.99, 0.75) rectangle (2.99, 1.49);
		\draw (0.5, 1.125) node {\includegraphics[width=7cm]{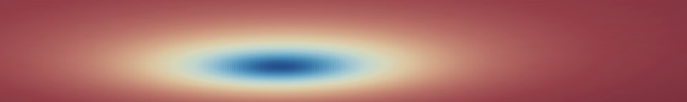}};

		% hemisphere shadow
		%\addplot[samples=45, domain=0:360, color = mygreylight, fill = mygreylight] (
		%	{cos(x)*0.5}, {0.15*sin(x)*0.5+1}
		%);
	
		% hemisphere lower silhouette
		\addplot[samples=30, domain=pi:2*pi, color = mygrey, densely dotted] (
			{cos(deg(x))}, {sin(deg(x))}
		);
		% hemisphere silhouette
		\addplot[samples=30, domain=0:pi, color = mygrey, fill = mygreylighter] (
			{cos(deg(x))}, {sin(deg(x))}
		);
		% hemisphere base
		\addplot[samples=25, domain=0:360, color = mygreylighterr, fill = mygreylighterr] (
			{cos(x)}, {0.15*sin(x)}
		);
		% projection helper
		\draw (axis cs:1.48,1.1) -- (axis cs:0,0) [mygreylight, dashed];
		\draw (axis cs:2.225,1) -- (axis cs:0,0) [mygreylight, dashed];

		% infinitesimal solid angle
		\addplot[
			rotate around={-60:(current axis.origin)},
			samples=45,
			domain=0:360,
			color = mygreen,
			fill = mygreenlight
		] (
			{0.4*0.25*cos(x)}, {0.4*0.125*sin(x)+0.85}
		);
		% infinitesimal area
		\addplot[
			rotate around={2:(current axis.origin)},
			samples=45,
			domain=0:360,
			color = myredlight
		] (
			{0.4*cos(x)+1.94}, {0.125*sin(x)+1.025}
		);
		% hemisphere base
		\addplot[samples=45, domain=180:360, color = mygreylight, densely dotted] (
			{cos(x)}, {0.15*sin(x)}
		);
		\addplot[samples=45, domain=0:180, color = mygrey] (
			{cos(x)}, {0.15*sin(x)}
		);
		% projection center
		\draw[mygrey, fill=mygreylighter] (axis cs:0,0) circle [radius=0.125*0.4];
		% labels
		%\node[black, draw = none, anchor = center] at (axis cs:1.35, 0.05)
		%	{\scalebox{0.5}{$\theta$}};
		\node[myredlight, draw = none, anchor = center] at (axis cs:2.0, 1.325)
			{\scalebox{0.7}{$d \bx$}};
		\draw[mygreen, ->, > = stealth] (axis cs: 1.08, 0.20) to [bend left = 10] (axis cs: 0.8, 0.33);
		\node[mygreen, draw = none, anchor = north west] at (axis cs:1.0, 0.38)
			{\scalebox{0.7}{$d \bomega$}};
		% Sets
		\node[mygrey, draw = none, anchor = center] at (axis cs:-0.6, -0.55)
			{\scalebox{0.8}{$\mathcal S ^2$}};
		\node[mygreylight, draw = none, anchor = center, xslant=0.0] at (axis cs:-1.55, 1.15)
			{\scalebox{0.8}{$\mathbb R^2 \times \{1\}$}};
		% formulas
		\draw[xshift = 5.5, yshift = -2.0, rounded corners=1, mygrey, fill = mygreylighter, drop shadow] 
		(0.85, -1.1) rectangle (2.5, -0.5);
		\node[xshift = 4.0, yshift = -2.0, black, draw = none, anchor = north west] at (0.80, -0.60)
			{\scalebox{0.7}{
				${\color{mygreen} d\bomega} =  \frac{1}{(x_1^2+x_2^2+1)^{3/2}} {\color{myred} d \bx}$
			}};
    \end{axis}
    \begin{pgfonlayer}{background}
    \draw[fill=white, drop shadow]
        (myaxis.north west) 
        rectangle (myaxis.south east);
    \end{pgfonlayer}
\end{tikzpicture}

%% file: fig_pdf_std.tex
\pgfplotsset{every axis/.append style = {
	x = 0.7cm,
	y = 0.7cm,
	xmin = -5, xmax = 5,
	ymin = -3, ymax = 3,
	grid = both,
	axis line style = {color = black},
	minor grid style={densely dotted},
	major grid style={densely dotted},
	y label style = {rotate = -90},
	minor tick num = 0,
	xtick = {-5,...,5},
	ytick = {-3,...,3},
	xticklabels = {},
	yticklabels = {},
	no markers
}}
\begin{tikzpicture}
    \begin{axis}[name = myaxis, axis on top] 
		\addplot graphics[xmin=-5,ymin=-3,xmax=5,ymax=3] {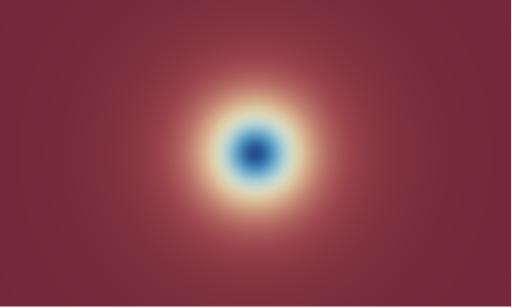};
	\end{axis}
    \begin{pgfonlayer}{background}
    \draw[fill=white, drop shadow]
        (myaxis.north west) 
        rectangle (myaxis.south east);
    \end{pgfonlayer}
\end{tikzpicture}

%% file: fig_pdf_int.tex
\pgfplotsset{every axis/.append style = {
	x = 0.7cm,
	y = 0.7cm,
	xmin = -5, xmax = 5,
	ymin = -3, ymax = 3,
	grid = both,
	axis line style = {color = black},
	minor grid style={densely dotted},
	major grid style={densely dotted},
	y label style = {rotate = -90},
	minor tick num = 0,
	xtick = {-5,...,5},
	ytick = {-3,...,3},
	xticklabels = {},
	yticklabels = {},
	no markers
}}
\begin{tikzpicture}
	\begin{axis}[name = myaxis, axis on top] 
		\coordinate (A) at (-2.7, 1.4);
		\coordinate (B) at (-1.2, -1.2);
		\coordinate (C) at (3.0, 2.3);
		\coordinate (D) at (-1.76, 1.26);
		\addplot graphics[xmin=-5,ymin=-3,xmax=5,ymax=3] {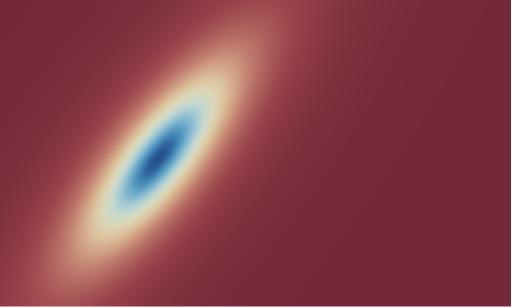};
		\filldraw[white, opacity = 0.8, even odd rule]	
			(axis cs: -5, -3) rectangle (axis cs: 5, 3)
			(A) -- (B) -- (C) -- cycle;
		\draw (A) -- (B) -- (C) -- cycle;
		%\draw[black, fill=mygreylight] (A) circle [radius=0.125*0.2];
		%\draw[black, fill=mygreylight] (B) circle [radius=0.125*0.2];
		%\draw[black, fill=mygreylight] (C) circle [radius=0.125*0.2];
		%\draw[black, fill=mygreylight] (D) circle [radius=0.125*0.2];
		\node[black,draw=none] at (2.5, 1.5) {\scalebox{0.75}{$\mathbb{D}$}};

		%samples 
		\draw[black, fill=mygreylight] (-2.59781, 1.32601) circle [radius=0.125*0.2];
		\draw[black, fill=mygreylight] (0.593813, 1.87004) circle [radius=0.125*0.2];
		\draw[black, fill=mygreylight] (0.637932, 0.368922) circle [radius=0.125*0.2];
		\draw[black, fill=mygreylight] (-0.302005, 1.67593) circle [radius=0.125*0.2];
		\draw[black, fill=mygreylight] (-0.355196, 1.69385) circle [radius=0.125*0.2];
		\draw[black, fill=mygreylight] (-0.222058, -0.111651) circle [radius=0.125*0.2];
		\draw[black, fill=mygreylight] (-0.688983, 1.62578) circle [radius=0.125*0.2];
		\draw[black, fill=mygreylight] (-0.809413, 1.67792) circle [radius=0.125*0.2];
		\draw[black, fill=mygreylight] (-0.889716, -0.409783) circle [radius=0.125*0.2];
		\draw[black, fill=mygreylight] (-0.996195, 1.62108) circle [radius=0.125*0.2];
		\draw[black, fill=mygreylight] (-1.3108, 1.5976) circle [radius=0.125*0.2];
		\draw[black, fill=mygreylight] (-1.01661, -0.715167) circle [radius=0.125*0.2];
		\draw[black, fill=mygreylight] (-0.962935, -0.975764) circle [radius=0.125*0.2];
		\draw[black, fill=mygreylight] (-1.70259, 1.39688) circle [radius=0.125*0.2];
		\draw[black, fill=mygreylight] (-1.18052, -0.985375) circle [radius=0.125*0.2];
		\draw[black, fill=mygreylight] (-2.55372, 1.30963) circle [radius=0.125*0.2];
		\draw[black, fill=mygreylight] (1.28763, 1.10692) circle [radius=0.125*0.2];
		\draw[black, fill=mygreylight] (-1.39543, -0.807791) circle [radius=0.125*0.2];
		\draw[black, fill=mygreylight] (1.43245, 1.08837) circle [radius=0.125*0.2];
		\draw[black, fill=mygreylight] (0.177406, 0.589423) circle [radius=0.125*0.2];
		\draw[black, fill=mygreylight] (-0.288377, 0.217474) circle [radius=0.125*0.2];
		\draw[black, fill=mygreylight] (-0.421871, 0.13209) circle [radius=0.125*0.2];
		\draw[black, fill=mygreylight] (-0.43979, -0.0343746) circle [radius=0.125*0.2];
		\draw[black, fill=mygreylight] (-0.289442, -0.367039) circle [radius=0.125*0.2];
		\draw[black, fill=mygreylight] (-0.843085, -0.220346) circle [radius=0.125*0.2];
		\draw[black, fill=mygreylight] (-0.987673, -0.206761) circle [radius=0.125*0.2];
		\draw[black, fill=mygreylight] (-1.06079, -0.326842) circle [radius=0.125*0.2];
		\draw[black, fill=mygreylight] (-1.04293, -0.535176) circle [radius=0.125*0.2];
		\draw[black, fill=mygreylight] (-0.825139, -0.760913) circle [radius=0.125*0.2];
		\draw[black, fill=mygreylight] (-1.27822, -0.496274) circle [radius=0.125*0.2];
		\draw[black, fill=mygreylight] (-1.04804, -0.844607) circle [radius=0.125*0.2];
		\draw[black, fill=mygreylight] (-1.48874, -0.556736) circle [radius=0.125*0.2];
		\draw[black, fill=mygreylight] (-1.5436, -0.581956) circle [radius=0.125*0.2];
		\draw[black, fill=mygreylight] (0.785254, 0.827909) circle [radius=0.125*0.2];
		\draw[black, fill=mygreylight] (0.293293, 0.804068) circle [radius=0.125*0.2];
		\draw[black, fill=mygreylight] (2.16454, 1.73551) circle [radius=0.125*0.2];
		\draw[black, fill=mygreylight] (0.178953, 0.498516) circle [radius=0.125*0.2];
		\draw[black, fill=mygreylight] (0.000934677, 0.301603) circle [radius=0.125*0.2];
		\draw[black, fill=mygreylight] (-0.501797, 0.0690341) circle [radius=0.125*0.2];
		\draw[black, fill=mygreylight] (-0.503615, 0.0780572) circle [radius=0.125*0.2];
		\draw[black, fill=mygreylight] (-0.552047, -0.0799762) circle [radius=0.125*0.2];
		\draw[black, fill=mygreylight] (-0.869765, -0.216345) circle [radius=0.125*0.2];
		\draw[black, fill=mygreylight] (-0.982946, -0.175899) circle [radius=0.125*0.2];
		\draw[black, fill=mygreylight] (-1.22958, -0.280122) circle [radius=0.125*0.2];
		\draw[black, fill=mygreylight] (-1.35388, -0.228133) circle [radius=0.125*0.2];
		\draw[black, fill=mygreylight] (-1.48842, -0.175085) circle [radius=0.125*0.2];
		\draw[black, fill=mygreylight] (-1.40618, -0.257376) circle [radius=0.125*0.2];
		\draw[black, fill=mygreylight] (-1.52991, -0.283825) circle [radius=0.125*0.2];
		\draw[black, fill=mygreylight] (-1.58591, -0.315707) circle [radius=0.125*0.2];
		\draw[black, fill=mygreylight] (-1.63305, -0.438171) circle [radius=0.125*0.2];
		\draw[black, fill=mygreylight] (-1.57721, -0.446401) circle [radius=0.125*0.2];
		\draw[black, fill=mygreylight] (-1.74007, -0.234834) circle [radius=0.125*0.2];
		\draw[black, fill=mygreylight] (0.150479, 0.823836) circle [radius=0.125*0.2];
		\draw[black, fill=mygreylight] (-0.358816, 0.289927) circle [radius=0.125*0.2];
		\draw[black, fill=mygreylight] (-0.674787, 0.387775) circle [radius=0.125*0.2];
		\draw[black, fill=mygreylight] (-0.783568, 0.202049) circle [radius=0.125*0.2];
		\draw[black, fill=mygreylight] (-0.886054, 0.138556) circle [radius=0.125*0.2];
		\draw[black, fill=mygreylight] (-0.874458, 0.00431988) circle [radius=0.125*0.2];
		\draw[black, fill=mygreylight] (-1.06561, 0.0731495) circle [radius=0.125*0.2];
		\draw[black, fill=mygreylight] (-1.25783, 0.0410157) circle [radius=0.125*0.2];
		\draw[black, fill=mygreylight] (-1.36973, 0.039418) circle [radius=0.125*0.2];
		\draw[black, fill=mygreylight] (-1.42997, -0.113589) circle [radius=0.125*0.2];
		\draw[black, fill=mygreylight] (-1.54836, -0.0666788) circle [radius=0.125*0.2];
		\draw[black, fill=mygreylight] (-1.42187, -0.316832) circle [radius=0.125*0.2];
		\draw[black, fill=mygreylight] (-1.56783, -0.22079) circle [radius=0.125*0.2];
		\draw[black, fill=mygreylight] (-1.65503, -0.252425) circle [radius=0.125*0.2];
		\draw[black, fill=mygreylight] (-1.72079, -0.134907) circle [radius=0.125*0.2];
		\draw[black, fill=mygreylight] (-1.61371, -0.466994) circle [radius=0.125*0.2];
		\draw[black, fill=mygreylight] (0.876506, 1.22391) circle [radius=0.125*0.2];
		\draw[black, fill=mygreylight] (0.152811, 0.820848) circle [radius=0.125*0.2];
		\draw[black, fill=mygreylight] (-0.49157, 0.524057) circle [radius=0.125*0.2];
		\draw[black, fill=mygreylight] (-0.546541, 0.521233) circle [radius=0.125*0.2];
		\draw[black, fill=mygreylight] (-0.943134, 0.358344) circle [radius=0.125*0.2];
		\draw[black, fill=mygreylight] (-0.917014, 0.272581) circle [radius=0.125*0.2];
		\draw[black, fill=mygreylight] (-1.10106, 0.050914) circle [radius=0.125*0.2];
		\draw[black, fill=mygreylight] (-1.04425, 0.0724818) circle [radius=0.125*0.2];
		\draw[black, fill=mygreylight] (-1.32596, -0.0146021) circle [radius=0.125*0.2];
		\draw[black, fill=mygreylight] (-1.4221, 0.0846141) circle [radius=0.125*0.2];
		\draw[black, fill=mygreylight] (-1.53961, 0.111747) circle [radius=0.125*0.2];
		\draw[black, fill=mygreylight] (-1.48645, 0.0300591) circle [radius=0.125*0.2];
		\draw[black, fill=mygreylight] (-1.53935, -0.00672729) circle [radius=0.125*0.2];
		\draw[black, fill=mygreylight] (-1.59947, -0.187675) circle [radius=0.125*0.2];
		\draw[black, fill=mygreylight] (-1.67621, -0.154693) circle [radius=0.125*0.2];
		\draw[black, fill=mygreylight] (-1.66104, -0.26292) circle [radius=0.125*0.2];
		\draw[black, fill=mygreylight] (-1.73519, -0.175046) circle [radius=0.125*0.2];
		\draw[black, fill=mygreylight] (-1.78081, -0.133449) circle [radius=0.125*0.2];
		\draw[black, fill=mygreylight] (1.60986, 1.66396) circle [radius=0.125*0.2];
		\draw[black, fill=mygreylight] (-0.53574, 0.659605) circle [radius=0.125*0.2];
		\draw[black, fill=mygreylight] (-0.615704, 0.561103) circle [radius=0.125*0.2];
		\draw[black, fill=mygreylight] (-0.898997, 0.409574) circle [radius=0.125*0.2];
		\draw[black, fill=mygreylight] (-0.996483, 0.418742) circle [radius=0.125*0.2];
		\draw[black, fill=mygreylight] (-0.981536, 0.271573) circle [radius=0.125*0.2];
		\draw[black, fill=mygreylight] (-1.17646, 0.372303) circle [radius=0.125*0.2];
		\draw[black, fill=mygreylight] (-1.30474, 0.267219) circle [radius=0.125*0.2];
		\draw[black, fill=mygreylight] (-1.44818, 0.119342) circle [radius=0.125*0.2];
		\draw[black, fill=mygreylight] (-1.51937, 0.212157) circle [radius=0.125*0.2];
		\draw[black, fill=mygreylight] (-1.56163, 0.121005) circle [radius=0.125*0.2];
		\draw[black, fill=mygreylight] (-1.64202, 0.0574495) circle [radius=0.125*0.2];
		\draw[black, fill=mygreylight] (-1.72165, 0.110068) circle [radius=0.125*0.2];
		\draw[black, fill=mygreylight] (-1.72834, 0.0768632) circle [radius=0.125*0.2];
		\draw[black, fill=mygreylight] (-1.83673, -0.0231701) circle [radius=0.125*0.2];
		\draw[black, fill=mygreylight] (1.42454, 1.52742) circle [radius=0.125*0.2];
		\draw[black, fill=mygreylight] (-1.88619, 0.00488203) circle [radius=0.125*0.2];
		\draw[black, fill=mygreylight] (1.29884, 1.50439) circle [radius=0.125*0.2];
		\draw[black, fill=mygreylight] (-0.521424, 0.637864) circle [radius=0.125*0.2];
		\draw[black, fill=mygreylight] (-0.481231, 0.628887) circle [radius=0.125*0.2];
		\draw[black, fill=mygreylight] (-0.854564, 0.456795) circle [radius=0.125*0.2];
		\draw[black, fill=mygreylight] (-1.14021, 0.511331) circle [radius=0.125*0.2];
		\draw[black, fill=mygreylight] (-1.05338, 0.447485) circle [radius=0.125*0.2];
		\draw[black, fill=mygreylight] (-1.18653, 0.500451) circle [radius=0.125*0.2];
		\draw[black, fill=mygreylight] (-1.30981, 0.256811) circle [radius=0.125*0.2];
		\draw[black, fill=mygreylight] (-1.45789, 0.235691) circle [radius=0.125*0.2];
		\draw[black, fill=mygreylight] (-1.55929, 0.326469) circle [radius=0.125*0.2];
		\draw[black, fill=mygreylight] (-1.52784, 0.250715) circle [radius=0.125*0.2];
		\draw[black, fill=mygreylight] (-1.64807, 0.249414) circle [radius=0.125*0.2];
		\draw[black, fill=mygreylight] (-1.79651, 0.168479) circle [radius=0.125*0.2];
		\draw[black, fill=mygreylight] (-1.7674, 0.162249) circle [radius=0.125*0.2];
		\draw[black, fill=mygreylight] (-1.85907, -0.0154368) circle [radius=0.125*0.2];
		\draw[black, fill=mygreylight] (-1.90685, 0.0751584) circle [radius=0.125*0.2];
		\draw[black, fill=mygreylight] (0.412512, 1.2148) circle [radius=0.125*0.2];
		\draw[black, fill=mygreylight] (-0.246693, 0.956123) circle [radius=0.125*0.2];
		\draw[black, fill=mygreylight] (-0.427797, 0.872687) circle [radius=0.125*0.2];
		\draw[black, fill=mygreylight] (-0.882592, 0.703077) circle [radius=0.125*0.2];
		\draw[black, fill=mygreylight] (-0.827292, 0.622014) circle [radius=0.125*0.2];
		\draw[black, fill=mygreylight] (-1.08974, 0.518319) circle [radius=0.125*0.2];
		\draw[black, fill=mygreylight] (-1.20884, 0.42062) circle [radius=0.125*0.2];
		\draw[black, fill=mygreylight] (-1.30311, 0.422911) circle [radius=0.125*0.2];
		\draw[black, fill=mygreylight] (-1.32907, 0.416363) circle [radius=0.125*0.2];
		\draw[black, fill=mygreylight] (-1.424, 0.338581) circle [radius=0.125*0.2];
		\draw[black, fill=mygreylight] (-1.46686, 0.339709) circle [radius=0.125*0.2];
		\draw[black, fill=mygreylight] (-1.55618, 0.382722) circle [radius=0.125*0.2];
		\draw[black, fill=mygreylight] (-1.63288, 0.318762) circle [radius=0.125*0.2];
		\draw[black, fill=mygreylight] (-1.68631, 0.231886) circle [radius=0.125*0.2];
		\draw[black, fill=mygreylight] (-1.8748, 0.319656) circle [radius=0.125*0.2];
		\draw[black, fill=mygreylight] (-1.93213, 0.264702) circle [radius=0.125*0.2];
		\draw[black, fill=mygreylight] (0.893595, 1.34738) circle [radius=0.125*0.2];
		\draw[black, fill=mygreylight] (1.19456, 1.5816) circle [radius=0.125*0.2];
		\draw[black, fill=mygreylight] (0.0510479, 1.13161) circle [radius=0.125*0.2];
		\draw[black, fill=mygreylight] (-0.104987, 1.05808) circle [radius=0.125*0.2];
		\draw[black, fill=mygreylight] (-0.669869, 0.874425) circle [radius=0.125*0.2];
		\draw[black, fill=mygreylight] (-0.945215, 0.772222) circle [radius=0.125*0.2];
		\draw[black, fill=mygreylight] (-0.906321, 0.661587) circle [radius=0.125*0.2];
		\draw[black, fill=mygreylight] (-1.16819, 0.645367) circle [radius=0.125*0.2];
		\draw[black, fill=mygreylight] (-1.31831, 0.639234) circle [radius=0.125*0.2];
		\draw[black, fill=mygreylight] (-1.40817, 0.446553) circle [radius=0.125*0.2];
		\draw[black, fill=mygreylight] (-1.38147, 0.608786) circle [radius=0.125*0.2];
		\draw[black, fill=mygreylight] (-1.52193, 0.549818) circle [radius=0.125*0.2];
		\draw[black, fill=mygreylight] (-1.66228, 0.494039) circle [radius=0.125*0.2];
		\draw[black, fill=mygreylight] (-1.78351, 0.438149) circle [radius=0.125*0.2];
		\draw[black, fill=mygreylight] (-1.7143, 0.350445) circle [radius=0.125*0.2];
		\draw[black, fill=mygreylight] (-1.81842, 0.412769) circle [radius=0.125*0.2];
		\draw[black, fill=mygreylight] (-2.01324, 0.230474) circle [radius=0.125*0.2];
		\draw[black, fill=mygreylight] (-1.98972, 0.20811) circle [radius=0.125*0.2];
		\draw[black, fill=mygreylight] (1.61684, 1.75771) circle [radius=0.125*0.2];
		\draw[black, fill=mygreylight] (0.760629, 1.48834) circle [radius=0.125*0.2];
		\draw[black, fill=mygreylight] (-0.116783, 1.18785) circle [radius=0.125*0.2];
		\draw[black, fill=mygreylight] (-0.604595, 1.00904) circle [radius=0.125*0.2];
		\draw[black, fill=mygreylight] (-0.852635, 0.904541) circle [radius=0.125*0.2];
		\draw[black, fill=mygreylight] (-0.972057, 0.717565) circle [radius=0.125*0.2];
		\draw[black, fill=mygreylight] (-1.06284, 0.678507) circle [radius=0.125*0.2];
		\draw[black, fill=mygreylight] (-1.289, 0.693572) circle [radius=0.125*0.2];
		\draw[black, fill=mygreylight] (-1.34206, 0.601763) circle [radius=0.125*0.2];
		\draw[black, fill=mygreylight] (-1.36327, 0.646253) circle [radius=0.125*0.2];
		\draw[black, fill=mygreylight] (-1.57205, 0.660427) circle [radius=0.125*0.2];
		\draw[black, fill=mygreylight] (-1.64238, 0.648359) circle [radius=0.125*0.2];
		\draw[black, fill=mygreylight] (-1.62586, 0.454727) circle [radius=0.125*0.2];
		\draw[black, fill=mygreylight] (-1.78748, 0.477285) circle [radius=0.125*0.2];
		\draw[black, fill=mygreylight] (-1.81766, 0.357954) circle [radius=0.125*0.2];
		\draw[black, fill=mygreylight] (-1.95026, 0.504464) circle [radius=0.125*0.2];
		\draw[black, fill=mygreylight] (0.59349, 1.33148) circle [radius=0.125*0.2];
		\draw[black, fill=mygreylight] (-2.04079, 0.487061) circle [radius=0.125*0.2];
		\draw[black, fill=mygreylight] (1.99562, 1.94311) circle [radius=0.125*0.2];
		\draw[black, fill=mygreylight] (-0.452012, 1.03532) circle [radius=0.125*0.2];
		\draw[black, fill=mygreylight] (-0.569631, 0.994587) circle [radius=0.125*0.2];
		\draw[black, fill=mygreylight] (-0.877151, 0.92028) circle [radius=0.125*0.2];
		\draw[black, fill=mygreylight] (-0.938157, 0.968497) circle [radius=0.125*0.2];
		\draw[black, fill=mygreylight] (-1.12811, 0.915171) circle [radius=0.125*0.2];
		\draw[black, fill=mygreylight] (-1.23714, 0.80041) circle [radius=0.125*0.2];
		\draw[black, fill=mygreylight] (-1.22851, 0.850698) circle [radius=0.125*0.2];
		\draw[black, fill=mygreylight] (-1.37791, 0.741456) circle [radius=0.125*0.2];
		\draw[black, fill=mygreylight] (-1.51186, 0.774504) circle [radius=0.125*0.2];
		\draw[black, fill=mygreylight] (-1.5303, 0.624472) circle [radius=0.125*0.2];
		\draw[black, fill=mygreylight] (-1.62987, 0.569869) circle [radius=0.125*0.2];
		\draw[black, fill=mygreylight] (-1.73147, 0.649659) circle [radius=0.125*0.2];
		\draw[black, fill=mygreylight] (-1.84302, 0.494756) circle [radius=0.125*0.2];
		\draw[black, fill=mygreylight] (-1.96952, 0.474997) circle [radius=0.125*0.2];
		\draw[black, fill=mygreylight] (0.286156, 1.36278) circle [radius=0.125*0.2];
		\draw[black, fill=mygreylight] (-2.24391, 0.690146) circle [radius=0.125*0.2];
		\draw[black, fill=mygreylight] (-0.0448337, 1.30507) circle [radius=0.125*0.2];
		\draw[black, fill=mygreylight] (-0.401711, 1.10623) circle [radius=0.125*0.2];
		\draw[black, fill=mygreylight] (-0.557548, 1.1083) circle [radius=0.125*0.2];
		\draw[black, fill=mygreylight] (-0.636892, 1.19563) circle [radius=0.125*0.2];
		\draw[black, fill=mygreylight] (-1.01078, 1.03874) circle [radius=0.125*0.2];
		\draw[black, fill=mygreylight] (-1.15598, 0.909013) circle [radius=0.125*0.2];
		\draw[black, fill=mygreylight] (-1.05696, 0.97195) circle [radius=0.125*0.2];
		\draw[black, fill=mygreylight] (-1.21187, 0.841542) circle [radius=0.125*0.2];
		\draw[black, fill=mygreylight] (-1.42755, 0.801753) circle [radius=0.125*0.2];
		\draw[black, fill=mygreylight] (-1.53316, 0.864036) circle [radius=0.125*0.2];
		\draw[black, fill=mygreylight] (-1.60491, 0.750865) circle [radius=0.125*0.2];
		\draw[black, fill=mygreylight] (-1.72369, 0.663448) circle [radius=0.125*0.2];
		\draw[black, fill=mygreylight] (-1.8471, 0.674185) circle [radius=0.125*0.2];
		\draw[black, fill=mygreylight] (-2.09232, 0.723565) circle [radius=0.125*0.2];
		\draw[black, fill=mygreylight] (-1.99223, 0.700206) circle [radius=0.125*0.2];
		\draw[black, fill=mygreylight] (-2.17101, 0.503611) circle [radius=0.125*0.2];
		\draw[black, fill=mygreylight] (0.216379, 1.4014) circle [radius=0.125*0.2];
		\draw[black, fill=mygreylight] (-0.159979, 1.3745) circle [radius=0.125*0.2];
		\draw[black, fill=mygreylight] (-0.222283, 1.36544) circle [radius=0.125*0.2];
		\draw[black, fill=mygreylight] (-0.474558, 1.19513) circle [radius=0.125*0.2];
		\draw[black, fill=mygreylight] (-0.810993, 1.18847) circle [radius=0.125*0.2];
		\draw[black, fill=mygreylight] (-0.774085, 1.14373) circle [radius=0.125*0.2];
		\draw[black, fill=mygreylight] (-1.02962, 1.18284) circle [radius=0.125*0.2];
		\draw[black, fill=mygreylight] (-1.11926, 1.01965) circle [radius=0.125*0.2];
		\draw[black, fill=mygreylight] (-1.22376, 1.11028) circle [radius=0.125*0.2];
		\draw[black, fill=mygreylight] (-1.36112, 1.03133) circle [radius=0.125*0.2];
		\draw[black, fill=mygreylight] (-1.32889, 1.08187) circle [radius=0.125*0.2];
		\draw[black, fill=mygreylight] (-1.58867, 0.876173) circle [radius=0.125*0.2];
		\draw[black, fill=mygreylight] (-1.73389, 0.747888) circle [radius=0.125*0.2];
		\draw[black, fill=mygreylight] (-1.89609, 0.906261) circle [radius=0.125*0.2];
		\draw[black, fill=mygreylight] (-2.11959, 0.788706) circle [radius=0.125*0.2];
		\draw[black, fill=mygreylight] (-2.03333, 0.77474) circle [radius=0.125*0.2];
		\draw[black, fill=mygreylight] (0.206451, 1.47315) circle [radius=0.125*0.2];
		\draw[black, fill=mygreylight] (0.53177, 1.68112) circle [radius=0.125*0.2];
		\draw[black, fill=mygreylight] (0.0990807, 1.58179) circle [radius=0.125*0.2];
		\draw[black, fill=mygreylight] (-0.33323, 1.40596) circle [radius=0.125*0.2];
		\draw[black, fill=mygreylight] (-0.576713, 1.41387) circle [radius=0.125*0.2];
		\draw[black, fill=mygreylight] (-0.78715, 1.16917) circle [radius=0.125*0.2];
		\draw[black, fill=mygreylight] (-0.570034, 1.4304) circle [radius=0.125*0.2];
		\draw[black, fill=mygreylight] (-0.877103, 1.21776) circle [radius=0.125*0.2];
		\draw[black, fill=mygreylight] (-1.0818, 1.21511) circle [radius=0.125*0.2];
		\draw[black, fill=mygreylight] (-1.174, 1.08298) circle [radius=0.125*0.2];
		\draw[black, fill=mygreylight] (-1.38537, 1.0287) circle [radius=0.125*0.2];
		\draw[black, fill=mygreylight] (-1.2971, 1.27043) circle [radius=0.125*0.2];
		\draw[black, fill=mygreylight] (-1.67318, 1.05126) circle [radius=0.125*0.2];
		\draw[black, fill=mygreylight] (-1.68467, 0.925684) circle [radius=0.125*0.2];
		\draw[black, fill=mygreylight] (-1.77574, 1.15862) circle [radius=0.125*0.2];
		\draw[black, fill=mygreylight] (-2.10517, 0.796121) circle [radius=0.125*0.2];
		\draw[black, fill=mygreylight] (-2.10944, 0.906835) circle [radius=0.125*0.2];
		\draw[black, fill=mygreylight] (-2.32667, 0.831966) circle [radius=0.125*0.2];
		\draw[black, fill=mygreylight] (-2.41177, 1.03276) circle [radius=0.125*0.2];
		\draw[black, fill=mygreylight] (1.82439, 1.98739) circle [radius=0.125*0.2];
		\draw[black, fill=mygreylight] (-0.180252, 1.46693) circle [radius=0.125*0.2];
		\draw[black, fill=mygreylight] (-0.267902, 1.47282) circle [radius=0.125*0.2];
		\draw[black, fill=mygreylight] (-0.402179, 1.58978) circle [radius=0.125*0.2];
		\draw[black, fill=mygreylight] (-0.766493, 1.4525) circle [radius=0.125*0.2];
		\draw[black, fill=mygreylight] (-1.02846, 1.26744) circle [radius=0.125*0.2];
		\draw[black, fill=mygreylight] (-1.06156, 1.30522) circle [radius=0.125*0.2];
		\draw[black, fill=mygreylight] (-1.1384, 1.23985) circle [radius=0.125*0.2];
		\draw[black, fill=mygreylight] (-1.15122, 1.2928) circle [radius=0.125*0.2];
		\draw[black, fill=mygreylight] (-1.26112, 1.40019) circle [radius=0.125*0.2];
		\draw[black, fill=mygreylight] (-1.53669, 1.33418) circle [radius=0.125*0.2];
		\draw[black, fill=mygreylight] (-1.70194, 1.28081) circle [radius=0.125*0.2];
		\draw[black, fill=mygreylight] (-1.7648, 1.16041) circle [radius=0.125*0.2];
		\draw[black, fill=mygreylight] (-2.06224, 1.06752) circle [radius=0.125*0.2];
		\draw[black, fill=mygreylight] (-2.40349, 1.11341) circle [radius=0.125*0.2];
		\draw[black, fill=mygreylight] (-2.4381, 1.09994) circle [radius=0.125*0.2];
		\draw[black, fill=mygreylight] (-2.61123, 1.3952) circle [radius=0.125*0.2];
		\draw[black, fill=mygreylight] (0.0532936, 1.71398) circle [radius=0.125*0.2];
		\draw[black, fill=mygreylight] (0.123162, 1.82944) circle [radius=0.125*0.2];
		\draw[black, fill=mygreylight] (-0.379566, 1.65395) circle [radius=0.125*0.2];
		\draw[black, fill=mygreylight] (-0.410405, 1.66583) circle [radius=0.125*0.2];
		\draw[black, fill=mygreylight] (-0.732417, 1.62873) circle [radius=0.125*0.2];
		\draw[black, fill=mygreylight] (-0.733117, 1.45677) circle [radius=0.125*0.2];
		\draw[black, fill=mygreylight] (-1.01154, 1.36645) circle [radius=0.125*0.2];
		\draw[black, fill=mygreylight] (-1.06801, 1.43874) circle [radius=0.125*0.2];
		\draw[black, fill=mygreylight] (-1.00656, 1.59564) circle [radius=0.125*0.2];
		\draw[black, fill=mygreylight] (-1.26953, 1.39229) circle [radius=0.125*0.2];
		\draw[black, fill=mygreylight] (-1.3326, 1.55661) circle [radius=0.125*0.2];
		\draw[black, fill=mygreylight] (-1.53862, 1.39623) circle [radius=0.125*0.2];
		\draw[black, fill=mygreylight] (-1.68396, 1.21838) circle [radius=0.125*0.2];
		\draw[black, fill=mygreylight] (-2.2032, 1.34773) circle [radius=0.125*0.2];
		\draw[black, fill=mygreylight] (-2.17758, 1.44488) circle [radius=0.125*0.2];
		\draw[black, fill=mygreylight] (-2.39518, 1.25213) circle [radius=0.125*0.2];
		\draw[black, fill=mygreylight] (0.588848, 1.88083) circle [radius=0.125*0.2];
		\draw[black, fill=mygreylight] (1.26591, 0.948813) circle [radius=0.125*0.2];
		\draw[black, fill=mygreylight] (1.07921, 0.699974) circle [radius=0.125*0.2];
		\draw[black, fill=mygreylight] (-0.148214, 1.79402) circle [radius=0.125*0.2];
		\draw[black, fill=mygreylight] (-0.248658, 1.76109) circle [radius=0.125*0.2];
		\draw[black, fill=mygreylight] (-0.651592, 1.64068) circle [radius=0.125*0.2];
		\draw[black, fill=mygreylight] (-0.6339, 1.71516) circle [radius=0.125*0.2];
		\draw[black, fill=mygreylight] (-0.877828, 1.53228) circle [radius=0.125*0.2];
		\draw[black, fill=mygreylight] (-1.07654, 1.60862) circle [radius=0.125*0.2];
		\draw[black, fill=mygreylight] (-1.12072, 1.63374) circle [radius=0.125*0.2];
		\draw[black, fill=mygreylight] (-1.18033, 1.60671) circle [radius=0.125*0.2];
		\draw[black, fill=mygreylight] (-1.39484, 1.60265) circle [radius=0.125*0.2];
		\draw[black, fill=mygreylight] (-1.56845, 1.38834) circle [radius=0.125*0.2];
		\draw[black, fill=mygreylight] (-0.9357, -0.950002) circle [radius=0.125*0.2];
		\draw[black, fill=mygreylight] (-1.09783, -0.86231) circle [radius=0.125*0.2];
		\draw[black, fill=mygreylight] (-1.27346, -0.980734) circle [radius=0.125*0.2];

	\end{axis}
    \begin{pgfonlayer}{background}
    \draw[fill=white, drop shadow]
        (myaxis.north west) 
        rectangle (myaxis.south east);
    \end{pgfonlayer}
\end{tikzpicture}

%% file: fig_ctor_gnomonic_2d_infinitesimals.tex
\pgfplotsset{every axis/.append style = {
	x = 1.4cm,
	y = 1.4cm,
	xmin = -2, xmax = 3,
	ymin = -1.5, ymax = 1.5,
	grid = both,
	axis line style = {color = black},
	minor grid style={densely dotted},
	major grid style={densely dotted},
	y label style = {rotate = -90},
	minor tick num = 1,
	xtick = {-3,...,3},
	ytick = {-2,...,2},
	xticklabels = {},
	yticklabels = {},
	no markers
}}
\begin{tikzpicture}
    \begin{axis}[name = myaxis]
		% real plane
		\draw[mygreylight, fill = mygreylighterr] 
			(-1.98, 0.75) rectangle (2.98, 1.49);

		% hemisphere shadow
		\addplot[samples=45, domain=0:360, color = mygreylighter, fill = mygreylighter] (
			{cos(x)*0.5}, {0.15*sin(x)*0.5+1}
		);
		%\draw[opacity = 0.5] (0.5, 1.125) node {\includegraphics[width=7cm]{./densityPerspective.pdf}};
				
		% hemisphere lower silhouette
		\addplot[samples=30, domain=pi:2*pi, color = mygrey, densely dotted] (
			{cos(deg(x))}, {sin(deg(x))}
		);
		% hemisphere silhouette
		\addplot[samples=30, domain=0:pi, color = mygrey, fill = mygreylighter] (
			{cos(deg(x))}, {sin(deg(x))}
		);
		% hemisphere base
		\addplot[samples=25, domain=0:360, color = mygreylighterr, fill = mygreylighterr] (
			{cos(x)}, {0.15*sin(x)}
		);
		% projection helper
		\draw (axis cs:1.48,1.1) -- (axis cs:0,0) [mygreylight, dashed];
		\draw (axis cs:2.225,1) -- (axis cs:0,0) [mygreylight, dashed];

		% infinitesimal solid angle
		\addplot[
			rotate around={-60:(current axis.origin)},
			samples=45,
			domain=0:360,
			color = mygreenlight,
			fill = mygreenlight
		] (
			{0.4*0.25*cos(x)}, {0.4*0.125*sin(x)+0.85}
		);
		% infinitesimal area
		\addplot[
			rotate around={2:(current axis.origin)},
			samples=45,
			domain=0:360,
			color = myredlight,
			fill = myredlight
		] (
			{0.4*cos(x)+1.94}, {0.125*sin(x)+1.025}
		);
		
		% line
		\draw (axis cs:0,0) -- (axis cs:1.9*0.455488*0.285,1.1*0.455488*0.285) [black];
		\draw (axis cs:1.9*0.455488*0.85,1.1*0.455488*0.85) 
		-- (axis cs:1.9*0.455488*0.29,1.1*0.455488*0.29) [black, opacity = 0.5];
		% point
		\draw[myred, fill=myred] (1.9,1.1) circle [radius=0.125*0.2];
		\draw[mygreen, fill=mygreen] (1.9*0.455488*0.85,1.1*0.455488*0.85) circle [radius=0.125*0.2];
		
		% line ctd
		\draw (axis cs:1.9,1.1) 
		-- (axis cs:1.9*0.455488*0.85,1.1*0.455488*0.85) [black];

		% hemisphere base
		\addplot[samples=45, domain=180:360, color = mygreylight, densely dotted] (
			{cos(x)}, {0.15*sin(x)}
		);
		\addplot[samples=45, domain=0:180, color = mygrey] (
			{cos(x)}, {0.15*sin(x)}
		);
		% projection center
		\draw[mygrey, fill=mygreylighter] (axis cs:0,0) circle [radius=0.125*0.4];
		% labels
		%\node[black, draw = none, anchor = center] at (axis cs:1.35, 0.05)
		%	{\scalebox{0.5}{$\theta$}};
		\node[myredlight, draw = none, anchor = center] at (axis cs:2.0, 1.325)
			{\scalebox{0.7}{$d \bx$}};
		\node[myred, draw = none, anchor = center] at (axis cs:2.0, 1.1)
			{\scalebox{0.7}{$\bx$}};
		%\draw[mygreen, ->, > = stealth] (axis cs: 1.08, 0.20) to [bend left = 10] (axis cs: 0.8, 0.33);
		\node[mygreenlight, draw = none, anchor = north west] at (axis cs:0.37, 0.8)
			{\scalebox{0.7}{$d \bomega$}};
		\node[mygreen, draw = none, anchor = north west] at (0.65, 0.52)
			{\scalebox{0.7}{$\bomega$}};
		% Sets
		\node[mygreylight, draw = none, anchor = center] at (axis cs:-0.6, -0.55)
			{\scalebox{0.8}{$\mathcal S ^2$}};
		\node[mygreylight, draw = none, anchor = center, xslant=0.0] at (axis cs:-1.55, 1.15)
			{\scalebox{0.8}{$\mathbb R^2 \times \{1\}$}};
		% formulas
		\draw[xshift = 4.0, yshift = -2.0, rounded corners=1, mygrey, fill = mygreylighter, drop shadow] 
		(1, -1.1) rectangle (2.45, -0.5);
		\node[black, draw = none, anchor = north west] at (1.35, -0.52)
			{\scalebox{0.7}{
				${\color{mygreen} \bomega} =  g_\perp( {\color{myred}\bx} )$
			}};
		\node[black, draw = none, anchor = north west] at (1.15, -0.78)
			{\scalebox{0.7}{
				$\Rightarrow{\color{myred} \bx} =  g_\perp^{-1}( {\color{mygreen}\bomega} )$
			}};
    \end{axis}
    \begin{pgfonlayer}{background}
    \draw[fill=white, drop shadow]
        (myaxis.north west) 
        rectangle (myaxis.south east);
    \end{pgfonlayer}
\end{tikzpicture}

%% file: fig_ctor_gnomonic_2d_triangle.tex
\pgfplotsset{every axis/.append style = {
	x = 1.4cm,
	y = 1.4cm,
	xmin = -2, xmax = 3,
	ymin = -1.5, ymax = 1.5,
	grid = both,
	axis line style = {color = black},
	minor grid style={densely dotted},
	major grid style={densely dotted},
	y label style = {rotate = -90},
	minor tick num = 1,
	xtick = {-3,...,3},
	ytick = {-2,...,2},
	xticklabels = {},
	yticklabels = {},
	no markers
}}
\begin{tikzpicture}
    \begin{axis}[name = myaxis]
		% real plane
		\draw[mygreylight, fill = mygreylighterr] 
			(-1.99, 0.75) rectangle (2.99, 1.49);

		% hemisphere shadow
		\addplot[samples=45, domain=0:360, color = mygreylighter, fill = mygreylighter] (
			{cos(x)*0.5}, {0.15*sin(x)*0.5+1}
		);

		% density
		%\draw[opacity = 0.5] (0.5, 1.125) node {\includegraphics[width=7cm]{./densityPerspective.pdf}};
	
		% hemisphere lower silhouette
		\addplot[samples=30, domain=pi:2*pi, color = mygreylight, densely dotted] (
			{cos(deg(x))}, {sin(deg(x))}
		);
		% hemisphere silhouette
		\addplot[samples=30, domain=0:pi, color = mygreylight, fill = mygreylighterr] (
			{cos(deg(x))}, {sin(deg(x))}
		);
		% hemisphere base
		\addplot[samples=25, domain=0:360, color = mygreylighterr, fill = mygreylighterr] (
			{cos(x)}, {0.15*sin(x)}
		);
		% projection helper
		%\draw (axis cs:1.45,1) -- (axis cs:0,0) [mygreylight, dashed];
		%\draw (axis cs:2.25,1) -- (axis cs:0,0) [mygreylight, dashed];

		% infinitesimal solid angle
		\addplot[
			rotate around={-60:(current axis.origin)},
			samples=45,
			domain=0:360,
			color = mygreenlighter,
			fill = mygreenlighter
		] (
			{0.4*0.25*cos(x)}, {0.4*0.125*sin(x)+0.85}
		);
		% infinitesimal area
		\addplot[
			rotate around={2:(current axis.origin)},
			samples=45,
			domain=0:360,
			color = myredlighter,
			fill = myredlighter
		] (
			{0.4*cos(x)+1.94}, {0.125*sin(x)+1.025}
		);

		% angle helper
		\addplot[
			samples=12,
			domain=184:212.5,
			color = black,
			fill = black,
			fill opacity = 0.2
		] (
			{0.5*cos(x)+1.9}, {0.45*sin(x)+1.1}
		) -- (1.9, 1.1);

		% triangle
		\draw (axis cs:0,0) -- (axis cs:0,0.14) [black, line width=0.75];
		\draw (axis cs:0,0.16) -- (axis cs:0,1) [black, opacity = 0.5, line width=0.75];
		\draw (axis cs:0,0) -- (axis cs:1.9*0.455488*0.285,1.1*0.455488*0.285) [black, line width=0.75];
		\draw (axis cs:1.9*0.455488*0.85,1.1*0.455488*0.85) 
		-- (axis cs:1.9*0.455488*0.29,1.1*0.455488*0.29) [black, opacity = 0.5, line width=0.75];
		% point
		\draw[myredlight, fill=myredlight] (1.9,1.1) circle [radius=0.125*0.2];
		\draw[mygreenlight, fill=mygreenlight] (1.9*0.455488*0.85,1.1*0.455488*0.85) circle [radius=0.125*0.2];
		% triangle ctd
		\draw (axis cs:0,1) -- (axis cs:1.9,1.1) [black, line width=0.75];
		\draw (axis cs:1.9,1.1) 
		-- (axis cs:1.9*0.455488*0.85,1.1*0.455488*0.85) [black, line width=0.75];
		% perp helper
		\draw[black] (axis cs:0.05,1) -- (axis cs:0.05,1-0.07) -- (axis cs:0.0,1-0.075);		
		% hemisphere base
		\addplot[samples=45, domain=180:360, color = mygreylight, densely dotted] (
			{cos(x)}, {0.15*sin(x)}
		);
		\addplot[samples=45, domain=0:180, color = mygreylight] (
			{cos(x)}, {0.15*sin(x)}
		);
		
		% projection center
		\draw[mygrey, fill=mygreylighter] (axis cs:0,0) circle [radius=0.125*0.4];
		% labels
		\node[myredlighter, draw = none, anchor = center] at (axis cs:2.0, 1.325)
			{\scalebox{0.7}{$d \bx$}};
		\node[myredlight, draw = none, anchor = center] at (axis cs:2.0, 1.1)
			{\scalebox{0.7}{$\bx$}};
		%\draw[mygreenlight, ->, > = stealth] (axis cs: 1.08, 0.20) to [bend left = 10] (axis cs: 0.8, 0.33);
		\node[mygreenlighter, draw = none, anchor = north west] at (axis cs:0.37, 0.8)
			{\scalebox{0.7}{$d \bomega$}};
		\node[mygreenlight, draw = none, anchor = north west] at (0.65, 0.52)
			{\scalebox{0.7}{$\bomega$}};
		\node[black, draw = none, anchor = north west] at (axis cs:1.2, 1.08)
			{\scalebox{0.6}{$\theta$}};

		% measurement helper
		%\draw[black] 	(0.0 + 1.1*0.455488*0.07, 0.0 - 1.9*0.455488*0.07) 
		%		-- 		(0.0 + 1.1*0.455488*0.30, 0.0 - 1.9*0.455488*0.30);
		%\draw[black] 	(1.9 + 1.1*0.455488*0.02, 1.1 - 1.9*0.455488*0.02)
		%		-- 		(1.9 + 1.1*0.455488*0.30, 1.1 - 1.9*0.455488*0.30);
		\draw[black, ->|, > = stealth] (1.0, 0.4) -- (0.07, -0.15);
		\draw[black, ->|, > = stealth] (1.15, 0.5) -- (1.98, 0.99);

		\node[black, draw = none, anchor = west] at (0.91, 0.45)
			{\scalebox{0.6}{$\|\bx\| \,=\, \sqrt{x_1^2+x_2^2+1}$}};

		% Sets
		\node[mygreylight, draw = none, anchor = center] at (axis cs:-0.6, -0.55)
			{\scalebox{0.8}{$\mathcal S ^2$}};
		\node[mygreylight, draw = none, anchor = center, xslant=0.0] at (axis cs:-1.55, 1.15)
			{\scalebox{0.8}{$\mathbb R^2 \times \{1\}$}};
    \end{axis}
    \begin{pgfonlayer}{background}
    \draw[fill=white, drop shadow]
        (myaxis.north west) 
        rectangle (myaxis.south east);
    \end{pgfonlayer}
\end{tikzpicture}

%% file: fig_simulate.tex
\pgfplotsset{every axis/.append style = {
	x = 1.4cm,
	y = 1.4cm,
	xmin = -2, xmax = 3,
	ymin = -1.5, ymax = 1.5,
	grid = both,
	axis line style = {color = black},
	minor grid style={densely dotted},
	major grid style={densely dotted},
	y label style = {rotate = -90},
	minor tick num = 1,
	xtick = {-3,...,3},
	ytick = {-2,...,2},
	xticklabels = {},
	yticklabels = {},
	no markers
}}
\begin{tikzpicture}
    \begin{axis}[name = myaxis]
		% real plane
		\draw[white, fill = white] 
			(-2.99, 0.75) rectangle (3.99, 1.49);

		% density
		\draw[opacity = 1.0] (0.5, 1.125) node {\includegraphics[width=7cm]{./densityPerspective.pdf}};

		% Sets
		\node[mygreylight, draw = none, anchor = center] at (axis cs:-0.6, -0.55)
			{\scalebox{0.8}{$\mathcal S ^2$}};
		\node[black, draw = none, anchor = center, xslant=0.0] at (axis cs:-1.55, 1.15)
			{\scalebox{0.8}{$\mathbb R^2 \times \{1\}$}};

		% density mask
		\filldraw[white, opacity = 0.8, even odd rule]	
		(-3.985, 0.75) rectangle (3.985, 1.89)
		plot [smooth cycle] coordinates {
			(0.2, 1.2) (0.3, 1.32) (0.5, 1.4) (0.8, 1.44) (1.1, 1.42) (1.30, 1.37) (1.45, 1.3) (2.35, 1.37) 
			(2.42, 1.25) (2.0, 1.15) (1.1, 1.05) (0.6, 1.04) (0.25, 1.1)	
		};
		% mask outline
		\draw [black] plot [smooth cycle] coordinates {
			(0.2, 1.2) (0.3, 1.32) (0.5, 1.4) (0.8, 1.44) (1.1, 1.42) (1.30, 1.37) (1.45, 1.3) (2.35, 1.37) 
			(2.42, 1.25) (2.0, 1.15) (1.1, 1.05) (0.6, 1.04) (0.25, 1.1)	
		};

		% hemisphere shadow
		%\addplot[samples=45, domain=0:360, color = mygreylighter, fill = mygreylighter] (
		%	{cos(x)*0.5}, {0.15*sin(x)*0.5+1}
		%);
	
		% hemisphere lower silhouette
		\addplot[samples=30, domain=pi:2*pi, color = mygreylight, densely dotted] (
			{cos(deg(x))}, {sin(deg(x))}
		);
		
		% hemisphere silhouette
		\addplot[samples=30, domain=0:pi, color = mygreylight, fill = mygreylighter] (
			{cos(deg(x))}, {sin(deg(x))}
		);
		% hemisphere base
		\addplot[samples=25, domain=0:360, color = mygreylighterr, fill = mygreylighterr] (
			{cos(x)}, {0.15*sin(x)}
		);
		% projection helper
		\draw (0.18, 1.15) -- (0,0) [mygreylight, dashed];
		\draw (2.39, 1.23) -- (0,0) [mygreylight, dashed];
		% footprints
		\draw (0.5, 0) node {\includegraphics[width=7cm]{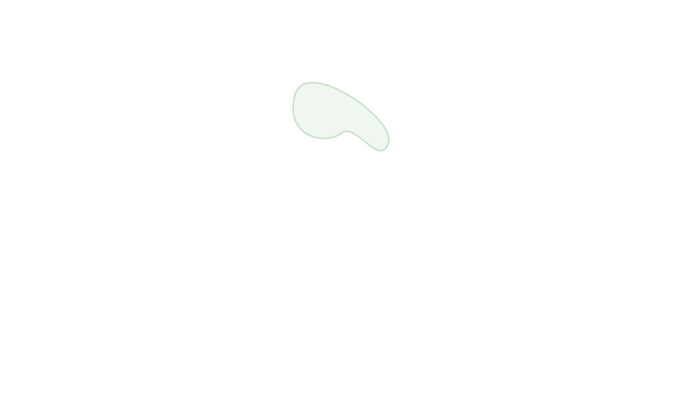}};

		% projective lines
		\draw (0, 0) -- (0.3,0.62) [mygrey];
		\draw (0, 0) -- (0.3,0.82) [mygrey];
		\draw (0, 0) -- (0.38,0.69) [mygrey];
		\draw (0, 0) -- (0.47,0.65) [mygrey];
		\draw (0, 0) -- (0.56,0.72) [mygrey];
		\draw (0, 0) -- (0.68,0.54) [mygrey];

		% hemisphere base
		\addplot[samples=45, domain=180:360, color = mygreylight, densely dotted] (
			{cos(x)}, {0.15*sin(x)}
		);
		\addplot[samples=45, domain=0:180, color = mygreylight] (
			{cos(x)}, {0.15*sin(x)}
		);

		% labels
		%\draw[mygreen, ->, > = stealth] (axis cs: 1.08, 0.20) to [bend left = 10] (axis cs: 0.8, 0.33);
		\draw[mygreen, ->, > = stealth] (1.08, 0.30) to [bend left = -10] (0.75, 0.46);
		\node[mygreenlight, draw = none, anchor = north west] at (-0.145, 1.0)
			{\scalebox{0.7}{$\Omega$}};
		\node[mygreen, draw = none, anchor = north west] at (1.21, 0.46)
			{\scalebox{0.7}{$\sim \mathcal U(\Omega) $}};
		\draw[mygreen, fill=mygreenlight] (1.2,0.28) circle [radius=0.125*0.2];
		\node[black, draw = none, anchor = north west] at (0.05, 1.55)
			{\scalebox{0.7}{$\mathbb D$}};
		\draw[black, ->, > = stealth] (1.98, 0.90) to [bend left = +10] (1.65, 1.20);
		\node[black, draw = none, anchor = north west] at (2.10, 1.05)
			{\scalebox{0.7}{$\sim \mathcal C(\mathbb D) $}};
		\draw[black, fill=mygreylight] (2.10,0.87) circle [radius=0.125*0.2];

		% random samples
		\draw[mygreen, fill=mygreenlight] (0.3,0.62) circle [radius=0.125*0.2];
		\draw[mygreen, fill=mygreenlight] (0.3,0.82) circle [radius=0.125*0.2];
		\draw[mygreen, fill=mygreenlight] (0.38,0.69) circle [radius=0.125*0.2];
		\draw[mygreen, fill=mygreenlight] (0.47,0.65) circle [radius=0.125*0.2];
		\draw[mygreen, fill=mygreenlight] (0.56,0.72) circle [radius=0.125*0.2];
		\draw[mygreen, fill=mygreenlight] (0.68,0.54) circle [radius=0.125*0.2];

		% random projected samples
		\draw[black, fill=mygreylight] (0.68,1.38) circle [radius=0.125*0.2];
		\draw[black, fill=mygreylight] (0.40,1.1) circle [radius=0.125*0.2];
		\draw[black, fill=mygreylight] (0.68,1.25) circle [radius=0.125*0.2];
		\draw[black, fill=mygreylight] (0.89,1.25) circle [radius=0.125*0.2];
		\draw[black, fill=mygreylight] (0.82,1.07) circle [radius=0.125*0.2];
		\draw[black, fill=mygreylight] (1.62,1.27) circle [radius=0.125*0.2];

		% projectives lines
		\draw (0.3,0.62) -- (0.68,1.38) [black];
		\draw (0.3,0.82)-- (0.4,1.1) [black];
		\draw (0.38,0.69) -- (0.68,1.25) [black];
		\draw (0.47,0.65) -- (0.89,1.25) [black];
		\draw (0.56,0.72) -- (0.82,1.07) [black];
		\draw (0.68,0.54) -- (1.62,1.27) [black];

		% projection center
		\draw[mygrey, fill=mygreylighter] (axis cs:0,0) circle [radius=0.125*0.4];

		% formulas
		\draw[xshift = 4.0, yshift = -2.0, rounded corners=1, mygrey, fill = mygreylighter, drop shadow] 
		(1, -1.1) rectangle (2.45, -0.5);
		\draw[mygreen, fill=mygreenlight] (1.35,-0.83) circle [radius=0.125*0.2];
		\draw[mygreen, fill=mygreenlight] (1.98,-0.83) circle [radius=0.125*0.2];
		\draw[black, fill=mygreylight] (2.35,-0.83) circle [radius=0.125*0.2];
		\node[xshift = 4.0, yshift = -2.0, black, draw = none, anchor = north west] at (1.25, -0.60)
			{\scalebox{0.7}{$ \mapsto g_\perp(\,\,) = $}};
    \end{axis}
    \begin{pgfonlayer}{background}
    \draw[fill=white, drop shadow]
        (myaxis.north west) 
        rectangle (myaxis.south east);
    \end{pgfonlayer}
\end{tikzpicture}

%% file: fig_integrate.tex
\pgfplotsset{every axis/.append style = {
	x = 1.4cm,
	y = 1.4cm,
	xmin = -2, xmax = 3,
	ymin = -1.5, ymax = 1.5,
	grid = both,
	axis line style = {color = black},
	minor grid style={densely dotted},
	major grid style={densely dotted},
	y label style = {rotate = -90},
	minor tick num = 1,
	xtick = {-3,...,3},
	ytick = {-2,...,2},
	xticklabels = {},
	yticklabels = {},
	no markers
}}
\begin{tikzpicture}
    \begin{axis}[name = myaxis]
		% real plane
		%\draw[mygreylight, fill = mygreylighterr] 
		%	(-1.99, 0.75) rectangle (2.99, 1.49);

		% density
		\draw (0.5, 1.125) node {\includegraphics[width=7cm]{./densityPerspective.pdf}};

		% Sets
		\node[mygreylight, draw = none, anchor = center] at (axis cs:-0.6, -0.55)
		{\scalebox{0.8}{$\mathcal S ^2$}};
		\node[black, draw = none, anchor = center, xslant=0.0] at (axis cs:-1.55, 1.15)
		{\scalebox{0.8}{$\mathbb R^2 \times \{1\}$}};

		% density mask
		\filldraw[white, opacity = 0.8, even odd rule]	
		(-2.985, 0.75) rectangle (3.985, 1.89)
		plot [smooth cycle] coordinates {
			(0.2, 1.2) (0.3, 1.32) (0.5, 1.4) (0.8, 1.44) (1.1, 1.42) (1.30, 1.37) (1.45, 1.3) (2.35, 1.37) 
			(2.42, 1.25) (2.0, 1.15) (1.1, 1.05) (0.6, 1.04) (0.25, 1.1)	
		};
		% mask outline
		\draw [black] plot [smooth cycle] coordinates {
			(0.2, 1.2) (0.3, 1.32) (0.5, 1.4) (0.8, 1.44) (1.1, 1.42) (1.30, 1.37) (1.45, 1.3) (2.35, 1.37) 
			(2.42, 1.25) (2.0, 1.15) (1.1, 1.05) (0.6, 1.04) (0.25, 1.1)	
		};

		% hemisphere shadow
		%\addplot[samples=45, domain=0:360, color = mygreylight, fill = mygreylight] (
		%	{cos(x)*0.5}, {0.15*sin(x)*0.5+1}
		%);
	
		% hemisphere lower silhouette
		\addplot[samples=30, domain=pi:2*pi, color = mygrey, densely dotted] (
			{cos(deg(x))}, {sin(deg(x))}
		);
		% hemisphere silhouette
		\addplot[samples=30, domain=0:pi, color = mygrey, fill = mygreylighter] (
			{cos(deg(x))}, {sin(deg(x))}
		);
		% hemisphere base
		\addplot[samples=25, domain=0:360, color = mygreylighterr, fill = mygreylighterr] (
			{cos(x)}, {0.15*sin(x)}
		);
		% projection helper
		\draw (0.18, 1.15) -- (0,0) [mygreylight, dashed];
		\draw (2.39, 1.23) -- (0,0) [mygreylight, dashed];
		% footprints
		\draw (0.5, 0) node {\includegraphics[width=7cm]{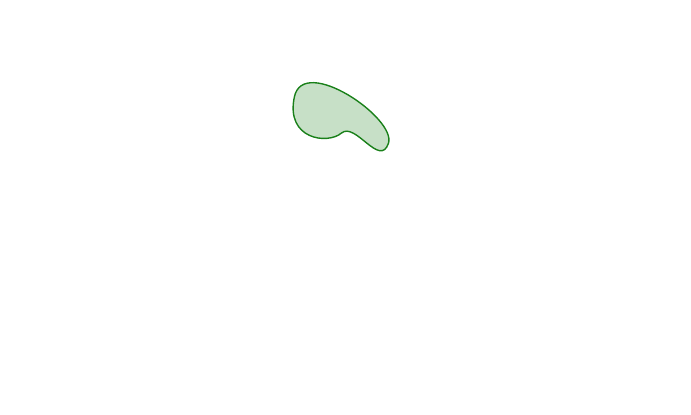}};

		% hemisphere base
		\addplot[samples=45, domain=180:360, color = mygreylight, densely dotted] (
			{cos(x)}, {0.15*sin(x)}
		);
		\addplot[samples=45, domain=0:180, color = mygrey] (
			{cos(x)}, {0.15*sin(x)}
		);
		% projection center
		\draw[mygrey, fill=mygreylighter] (axis cs:0,0) circle [radius=0.125*0.4];
		% labels
		%\node[black, draw = none, anchor = center] at (axis cs:1.35, 0.05)
		%	{\scalebox{0.5}{$\theta$}};
		%\draw[mygreen, ->, > = stealth] (axis cs: 1.08, 0.20) to [bend left = 10] (axis cs: 0.8, 0.33);
		\node[mygreen, draw = none, anchor = north west] at (-0.145, 1.0)
			{\scalebox{0.7}{$\Omega$}};
		\node[black, draw = none, anchor = north west] at (0.05, 1.55)
			{\scalebox{0.7}{$\mathbb D$}};

		% formulas
		\draw[xshift = 4.0, yshift = -2.0, rounded corners=1, mygrey, fill = mygreylighter, drop shadow] 
		(1, -1.1) rectangle (2.45, -0.5);
		\node[xshift = 4.0, yshift = -2.0, black, draw = none, anchor = north west] at (0.90, -0.60)
			{\scalebox{0.7}{
				${\color{black} \mathbb D} \mapsto g_\perp^{-1}({\color{black} \mathbb D}) = {\color{mygreen}\Omega}$
			}};
    \end{axis}
    \begin{pgfonlayer}{background}
    \draw[fill=white, drop shadow]
        (myaxis.north west) 
        rectangle (myaxis.south east);
    \end{pgfonlayer}
\end{tikzpicture}

%% file: fig_pdf.tex
\pgfplotsset{every axis/.append style = {
	x = 0.7cm,
	y = 0.7cm,
	xmin = -5, xmax = 5,
	ymin = -3, ymax = 3,
	grid = both,
	axis line style = {color = black},
	minor grid style={densely dotted},
	major grid style={densely dotted},
	y label style = {rotate = -90},
	minor tick num = 0,
	xtick = {-5,...,5},
	ytick = {-3,...,3},
	xticklabels = {},
	yticklabels = {},
	no markers
}}
\begin{tikzpicture}
    \begin{axis}[name = myaxis, axis on top] 
		\addplot graphics[xmin=-5,ymin=-3,xmax=5,ymax=3] {./pdf.jpg};
	\end{axis}
    \begin{pgfonlayer}{background}
    \draw[fill=white, drop shadow]
        (myaxis.north west) 
        rectangle (myaxis.south east);
    \end{pgfonlayer}
\end{tikzpicture}

%% file: fig_pdf_std_int.tex
\pgfplotsset{every axis/.append style = {
	x = 0.7cm,
	y = 0.7cm,
	xmin = -5, xmax = 5,
	ymin = -3, ymax = 3,
	grid = both,
	axis line style = {color = black},
	minor grid style={densely dotted},
	major grid style={densely dotted},
	y label style = {rotate = -90},
	minor tick num = 0,
	xtick = {-5,...,5},
	ytick = {-3,...,3},
	xticklabels = {},
	yticklabels = {},
	no markers
}}
\begin{tikzpicture}
    \begin{axis}[name = myaxis, axis on top] 
		\coordinate (A) at (-0.571429, 2.23249);
		\coordinate (B) at (0.5, -1.7451);
		\coordinate (C) at (3.5, -2.31373);
		\addplot graphics[xmin=-5,ymin=-3,xmax=5,ymax=3] {./pdf_std.jpg};
		\filldraw[white, opacity = 0.8, even odd rule]	
			(axis cs: -5, -3) rectangle (axis cs: 5, 3)
			(A) -- (B) -- (C) -- cycle;
		\draw (A) -- (B) -- (C) -- cycle;
		%\draw[black, fill=mygreylight] (A) circle [radius=0.125*0.2];
		%\draw[black, fill=mygreylight] (B) circle [radius=0.125*0.2];
		%\draw[black, fill=mygreylight] (C) circle [radius=0.125*0.2];
		%\draw[black, fill=mygreylight] (D) circle [radius=0.125*0.2];
		\node[black,draw=none] at (2.5, -0.5) {\scalebox{0.75}{$\mathbb{D}_\std$}};

		% samples
		\draw[black, fill=mygreylight] (-0.498435, 2.06263) circle [radius=0.125*0.2];
\draw[black, fill=mygreylight] (1.78129, -0.443651) circle [radius=0.125*0.2];
\draw[black, fill=mygreylight] (1.81281, -1.95735) circle [radius=0.125*0.2];
\draw[black, fill=mygreylight] (1.14142, 0.219213) circle [radius=0.125*0.2];
\draw[black, fill=mygreylight] (1.10343, 0.287434) circle [radius=0.125*0.2];
\draw[black, fill=mygreylight] (1.19853, -1.60946) circle [radius=0.125*0.2];
\draw[black, fill=mygreylight] (0.865012, 0.538591) circle [radius=0.125*0.2];
\draw[black, fill=mygreylight] (0.778991, 0.7044) circle [radius=0.125*0.2];
\draw[black, fill=mygreylight] (0.721631, -1.26588) circle [radius=0.125*0.2];
\draw[black, fill=mygreylight] (0.645575, 0.826571) circle [radius=0.125*0.2];
\draw[black, fill=mygreylight] (0.420858, 1.10317) circle [radius=0.125*0.2];
\draw[black, fill=mygreylight] (0.630996, -1.44443) circle [radius=0.125*0.2];
\draw[black, fill=mygreylight] (0.669332, -1.75103) circle [radius=0.125*0.2];
\draw[black, fill=mygreylight] (0.141005, 1.27953) circle [radius=0.125*0.2];
\draw[black, fill=mygreylight] (0.513917, -1.55324) circle [radius=0.125*0.2];
\draw[black, fill=mygreylight] (-0.466946, 2.00459) circle [radius=0.125*0.2];
\draw[black, fill=mygreylight] (2.27688, -1.85258) circle [radius=0.125*0.2];
\draw[black, fill=mygreylight] (0.360405, -1.17445) circle [radius=0.125*0.2];
\draw[black, fill=mygreylight] (2.38032, -2.0087) circle [radius=0.125*0.2];
\draw[black, fill=mygreylight] (1.48386, -1.30258) circle [radius=0.125*0.2];
\draw[black, fill=mygreylight] (1.15116, -1.22363) circle [radius=0.125*0.2];
\draw[black, fill=mygreylight] (1.05581, -1.1802) circle [radius=0.125*0.2];
\draw[black, fill=mygreylight] (1.04301, -1.32634) circle [radius=0.125*0.2];
\draw[black, fill=mygreylight] (1.1504, -1.79567) circle [radius=0.125*0.2];
\draw[black, fill=mygreylight] (0.754939, -1.12457) circle [radius=0.125*0.2];
\draw[black, fill=mygreylight] (0.651662, -0.973551) circle [radius=0.125*0.2];
\draw[black, fill=mygreylight] (0.599435, -1.02164) circle [radius=0.125*0.2];
\draw[black, fill=mygreylight] (0.612191, -1.2429) circle [radius=0.125*0.2];
\draw[black, fill=mygreylight] (0.767758, -1.67163) circle [radius=0.125*0.2];
\draw[black, fill=mygreylight] (0.44413, -0.980677) circle [radius=0.125*0.2];
\draw[black, fill=mygreylight] (0.608543, -1.5414) circle [radius=0.125*0.2];
\draw[black, fill=mygreylight] (0.29376, -0.839461) circle [radius=0.125*0.2];
\draw[black, fill=mygreylight] (0.254573, -0.811936) circle [radius=0.125*0.2];
\draw[black, fill=mygreylight] (1.91804, -1.64767) circle [radius=0.125*0.2];
\draw[black, fill=mygreylight] (1.56664, -1.20251) circle [radius=0.125*0.2];
\draw[black, fill=mygreylight] (2.90324, -2.07147) circle [radius=0.125*0.2];
\draw[black, fill=mygreylight] (1.48497, -1.39317) circle [radius=0.125*0.2];
\draw[black, fill=mygreylight] (1.35781, -1.41669) circle [radius=0.125*0.2];
\draw[black, fill=mygreylight] (0.998717, -1.1659) circle [radius=0.125*0.2];
\draw[black, fill=mygreylight] (0.997418, -1.15532) circle [radius=0.125*0.2];
\draw[black, fill=mygreylight] (0.962823, -1.26413) circle [radius=0.125*0.2];
\draw[black, fill=mygreylight] (0.735882, -1.09524) circle [radius=0.125*0.2];
\draw[black, fill=mygreylight] (0.655039, -0.947796) circle [radius=0.125*0.2];
\draw[black, fill=mygreylight] (0.478875, -0.81509) circle [radius=0.125*0.2];
\draw[black, fill=mygreylight] (0.390087, -0.645736) circle [radius=0.125*0.2];
\draw[black, fill=mygreylight] (0.293985, -0.465592) circle [radius=0.125*0.2];
\draw[black, fill=mygreylight] (0.352731, -0.624598) circle [radius=0.125*0.2];
\draw[black, fill=mygreylight] (0.26435, -0.532687) circle [radius=0.125*0.2];
\draw[black, fill=mygreylight] (0.224353, -0.510615) circle [radius=0.125*0.2];
\draw[black, fill=mygreylight] (0.19068, -0.585781) circle [radius=0.125*0.2];
\draw[black, fill=mygreylight] (0.230567, -0.647032) circle [radius=0.125*0.2];
\draw[black, fill=mygreylight] (0.114234, -0.284502) circle [radius=0.125*0.2];
\draw[black, fill=mygreylight] (1.46463, -1.04712) circle [radius=0.125*0.2];
\draw[black, fill=mygreylight] (1.10085, -1.08551) circle [radius=0.125*0.2];
\draw[black, fill=mygreylight] (0.875152, -0.688658) circle [radius=0.125*0.2];
\draw[black, fill=mygreylight] (0.797451, -0.767142) circle [radius=0.125*0.2];
\draw[black, fill=mygreylight] (0.724247, -0.731785) circle [radius=0.125*0.2];
\draw[black, fill=mygreylight] (0.73253, -0.874432) circle [radius=0.125*0.2];
\draw[black, fill=mygreylight] (0.595993, -0.624903) circle [radius=0.125*0.2];
\draw[black, fill=mygreylight] (0.458691, -0.473337) circle [radius=0.125*0.2];
\draw[black, fill=mygreylight] (0.378762, -0.368332) circle [radius=0.125*0.2];
\draw[black, fill=mygreylight] (0.335739, -0.460974) circle [radius=0.125*0.2];
\draw[black, fill=mygreylight] (0.251173, -0.302229) circle [radius=0.125*0.2];
\draw[black, fill=mygreylight] (0.341519, -0.667939) circle [radius=0.125*0.2];
\draw[black, fill=mygreylight] (0.237266, -0.434776) circle [radius=0.125*0.2];
\draw[black, fill=mygreylight] (0.174981, -0.382745) circle [radius=0.125*0.2];
\draw[black, fill=mygreylight] (0.128008, -0.204899) circle [radius=0.125*0.2];
\draw[black, fill=mygreylight] (0.204494, -0.632456) circle [radius=0.125*0.2];
\draw[black, fill=mygreylight] (1.98322, -1.34634) circle [radius=0.125*0.2];
\draw[black, fill=mygreylight] (1.46629, -1.05227) circle [radius=0.125*0.2];
\draw[black, fill=mygreylight] (1.00602, -0.729542) circle [radius=0.125*0.2];
\draw[black, fill=mygreylight] (0.966756, -0.679957) circle [radius=0.125*0.2];
\draw[black, fill=mygreylight] (0.683476, -0.461944) circle [radius=0.125*0.2];
\draw[black, fill=mygreylight] (0.702133, -0.570902) circle [radius=0.125*0.2];
\draw[black, fill=mygreylight] (0.570668, -0.612936) circle [radius=0.125*0.2];
\draw[black, fill=mygreylight] (0.611252, -0.645902) circle [radius=0.125*0.2];
\draw[black, fill=mygreylight] (0.410027, -0.462979) circle [radius=0.125*0.2];
\draw[black, fill=mygreylight] (0.341356, -0.274147) circle [radius=0.125*0.2];
\draw[black, fill=mygreylight] (0.257423, -0.135635) circle [radius=0.125*0.2];
\draw[black, fill=mygreylight] (0.295393, -0.266349) circle [radius=0.125*0.2];
\draw[black, fill=mygreylight] (0.25761, -0.252036) circle [radius=0.125*0.2];
\draw[black, fill=mygreylight] (0.214665, -0.372175) circle [radius=0.125*0.2];
\draw[black, fill=mygreylight] (0.159847, -0.26675) circle [radius=0.125*0.2];
\draw[black, fill=mygreylight] (0.170688, -0.387309) circle [radius=0.125*0.2];
\draw[black, fill=mygreylight] (0.117719, -0.230533) circle [radius=0.125*0.2];
\draw[black, fill=mygreylight] (0.0851326, -0.146303) circle [radius=0.125*0.2];
\draw[black, fill=mygreylight] (2.50705, -1.61336) circle [radius=0.125*0.2];
\draw[black, fill=mygreylight] (0.974471, -0.554585) circle [radius=0.125*0.2];
\draw[black, fill=mygreylight] (0.917354, -0.574999) circle [radius=0.125*0.2];
\draw[black, fill=mygreylight] (0.715002, -0.453754) circle [radius=0.125*0.2];
\draw[black, fill=mygreylight] (0.645369, -0.351922) circle [radius=0.125*0.2];
\draw[black, fill=mygreylight] (0.656046, -0.51044) circle [radius=0.125*0.2];
\draw[black, fill=mygreylight] (0.516817, -0.226047) circle [radius=0.125*0.2];
\draw[black, fill=mygreylight] (0.425184, -0.206893) circle [radius=0.125*0.2];
\draw[black, fill=mygreylight] (0.322725, -0.21526) circle [radius=0.125*0.2];
\draw[black, fill=mygreylight] (0.271878, -0.0564672) circle [radius=0.125*0.2];
\draw[black, fill=mygreylight] (0.241696, -0.10559) circle [radius=0.125*0.2];
\draw[black, fill=mygreylight] (0.184272, -0.0913334) circle [radius=0.125*0.2];
\draw[black, fill=mygreylight] (0.127395, 0.0360899) circle [radius=0.125*0.2];
\draw[black, fill=mygreylight] (0.122614, 0.00991007) circle [radius=0.125*0.2];
\draw[black, fill=mygreylight] (0.0451897, 0.0150705) circle [radius=0.125*0.2];
\draw[black, fill=mygreylight] (2.37467, -1.57072) circle [radius=0.125*0.2];
\draw[black, fill=mygreylight] (0.00986736, 0.089669) circle [radius=0.125*0.2];
\draw[black, fill=mygreylight] (2.28489, -1.47359) circle [radius=0.125*0.2];
\draw[black, fill=mygreylight] (0.984697, -0.589533) circle [radius=0.125*0.2];
\draw[black, fill=mygreylight] (1.01341, -0.636614) circle [radius=0.125*0.2];
\draw[black, fill=mygreylight] (0.74674, -0.449776) circle [radius=0.125*0.2];
\draw[black, fill=mygreylight] (0.54271, -0.124269) circle [radius=0.125*0.2];
\draw[black, fill=mygreylight] (0.604729, -0.269556) circle [radius=0.125*0.2];
\draw[black, fill=mygreylight] (0.509618, -0.0908136) circle [radius=0.125*0.2];
\draw[black, fill=mygreylight] (0.421563, -0.212269) circle [radius=0.125*0.2];
\draw[black, fill=mygreylight] (0.315793, -0.091948) circle [radius=0.125*0.2];
\draw[black, fill=mygreylight] (0.243362, 0.0936245) circle [radius=0.125*0.2];
\draw[black, fill=mygreylight] (0.265825, -0.0105951) circle [radius=0.125*0.2];
\draw[black, fill=mygreylight] (0.179953, 0.102625) circle [radius=0.125*0.2];
\draw[black, fill=mygreylight] (0.0739202, 0.164654) circle [radius=0.125*0.2];
\draw[black, fill=mygreylight] (0.0947111, 0.130825) circle [radius=0.125*0.2];
\draw[black, fill=mygreylight] (0.0292376, 0.0439216) circle [radius=0.125*0.2];
\draw[black, fill=mygreylight] (-0.00488943, 0.178243) circle [radius=0.125*0.2];
\draw[black, fill=mygreylight] (1.65179, -0.913371) circle [radius=0.125*0.2];
\draw[black, fill=mygreylight] (1.18093, -0.539164) circle [radius=0.125*0.2];
\draw[black, fill=mygreylight] (1.05157, -0.448484) circle [radius=0.125*0.2];
\draw[black, fill=mygreylight] (0.72672, -0.18163) circle [radius=0.125*0.2];
\draw[black, fill=mygreylight] (0.76622, -0.31377) circle [radius=0.125*0.2];
\draw[black, fill=mygreylight] (0.578759, -0.165483) circle [radius=0.125*0.2];
\draw[black, fill=mygreylight] (0.493687, -0.147837) circle [radius=0.125*0.2];
\draw[black, fill=mygreylight] (0.42635, -0.0558086) circle [radius=0.125*0.2];
\draw[black, fill=mygreylight] (0.407807, -0.0375042) circle [radius=0.125*0.2];
\draw[black, fill=mygreylight] (0.339997, -0.0233472) circle [radius=0.125*0.2];
\draw[black, fill=mygreylight] (0.309385, 0.0185744) circle [radius=0.125*0.2];
\draw[black, fill=mygreylight] (0.245588, 0.145806) circle [radius=0.125*0.2];
\draw[black, fill=mygreylight] (0.190803, 0.156148) circle [radius=0.125*0.2];
\draw[black, fill=mygreylight] (0.152638, 0.121861) circle [radius=0.125*0.2];
\draw[black, fill=mygreylight] (0.0179979, 0.38743) circle [radius=0.125*0.2];
\draw[black, fill=mygreylight] (-0.0229469, 0.388147) circle [radius=0.125*0.2];
\draw[black, fill=mygreylight] (1.99543, -1.24157) circle [radius=0.125*0.2];
\draw[black, fill=mygreylight] (2.2104, -1.29857) circle [radius=0.125*0.2];
\draw[black, fill=mygreylight] (1.39361, -0.650682) circle [radius=0.125*0.2];
\draw[black, fill=mygreylight] (1.28215, -0.574162) circle [radius=0.125*0.2];
\draw[black, fill=mygreylight] (0.878665, -0.216235) circle [radius=0.125*0.2];
\draw[black, fill=mygreylight] (0.681989, -0.0541991) circle [radius=0.125*0.2];
\draw[black, fill=mygreylight] (0.70977, -0.199706) circle [radius=0.125*0.2];
\draw[black, fill=mygreylight] (0.522719, 0.0337934) circle [radius=0.125*0.2];
\draw[black, fill=mygreylight] (0.415494, 0.170747) circle [radius=0.125*0.2];
\draw[black, fill=mygreylight] (0.351304, 0.0674309) circle [radius=0.125*0.2];
\draw[black, fill=mygreylight] (0.370377, 0.201052) circle [radius=0.125*0.2];
\draw[black, fill=mygreylight] (0.270052, 0.277006) circle [radius=0.125*0.2];
\draw[black, fill=mygreylight] (0.169802, 0.355989) circle [radius=0.125*0.2];
\draw[black, fill=mygreylight] (0.0832086, 0.416652) circle [radius=0.125*0.2];
\draw[black, fill=mygreylight] (0.132641, 0.264758) circle [radius=0.125*0.2];
\draw[black, fill=mygreylight] (0.0582691, 0.425022) circle [radius=0.125*0.2];
\draw[black, fill=mygreylight] (-0.0808856, 0.431842) circle [radius=0.125*0.2];
\draw[black, fill=mygreylight] (-0.0640874, 0.387518) circle [radius=0.125*0.2];
\draw[black, fill=mygreylight] (2.51203, -1.52809) circle [radius=0.125*0.2];
\draw[black, fill=mygreylight] (1.90045, -0.976739) circle [radius=0.125*0.2];
\draw[black, fill=mygreylight] (1.27373, -0.435707) circle [radius=0.125*0.2];
\draw[black, fill=mygreylight] (0.925289, -0.146421) circle [radius=0.125*0.2];
\draw[black, fill=mygreylight] (0.748118, -0.0126459) circle [radius=0.125*0.2];
\draw[black, fill=mygreylight] (0.662817, -0.0822214) circle [radius=0.125*0.2];
\draw[black, fill=mygreylight] (0.597969, -0.0340496) circle [radius=0.125*0.2];
\draw[black, fill=mygreylight] (0.436427, 0.196109) circle [radius=0.125*0.2];
\draw[black, fill=mygreylight] (0.398532, 0.156627) circle [radius=0.125*0.2];
\draw[black, fill=mygreylight] (0.383376, 0.220452) circle [radius=0.125*0.2];
\draw[black, fill=mygreylight] (0.234248, 0.433186) circle [radius=0.125*0.2];
\draw[black, fill=mygreylight] (0.184015, 0.488332) circle [radius=0.125*0.2];
\draw[black, fill=mygreylight] (0.195815, 0.282763) circle [radius=0.125*0.2];
\draw[black, fill=mygreylight] (0.0803731, 0.458802) circle [radius=0.125*0.2];
\draw[black, fill=mygreylight] (0.058817, 0.370552) circle [radius=0.125*0.2];
\draw[black, fill=mygreylight] (-0.0359033, 0.640483) circle [radius=0.125*0.2];
\draw[black, fill=mygreylight] (1.78106, -0.971336) circle [radius=0.125*0.2];
\draw[black, fill=mygreylight] (-0.100567, 0.709639) circle [radius=0.125*0.2];
\draw[black, fill=mygreylight] (2.78259, -1.70707) circle [radius=0.125*0.2];
\draw[black, fill=mygreylight] (1.03428, -0.265982) circle [radius=0.125*0.2];
\draw[black, fill=mygreylight] (0.950264, -0.193894) circle [radius=0.125*0.2];
\draw[black, fill=mygreylight] (0.730606, 0.0261326) circle [radius=0.125*0.2];
\draw[black, fill=mygreylight] (0.687031, 0.131506) circle [radius=0.125*0.2];
\draw[black, fill=mygreylight] (0.551353, 0.260128) circle [radius=0.125*0.2];
\draw[black, fill=mygreylight] (0.473474, 0.251456) circle [radius=0.125*0.2];
\draw[black, fill=mygreylight] (0.479638, 0.292539) circle [radius=0.125*0.2];
\draw[black, fill=mygreylight] (0.372923, 0.327726) circle [radius=0.125*0.2];
\draw[black, fill=mygreylight] (0.277239, 0.487704) circle [radius=0.125*0.2];
\draw[black, fill=mygreylight] (0.26407, 0.358174) circle [radius=0.125*0.2];
\draw[black, fill=mygreylight] (0.192951, 0.399467) circle [radius=0.125*0.2];
\draw[black, fill=mygreylight] (0.12038, 0.574453) circle [radius=0.125*0.2];
\draw[black, fill=mygreylight] (0.0407025, 0.528824) circle [radius=0.125*0.2];
\draw[black, fill=mygreylight] (-0.0496559, 0.62993) circle [radius=0.125*0.2];
\draw[black, fill=mygreylight] (1.56154, -0.647952) circle [radius=0.125*0.2];
\draw[black, fill=mygreylight] (-0.245652, 1.10219) circle [radius=0.125*0.2];
\draw[black, fill=mygreylight] (1.32512, -0.389306) circle [radius=0.125*0.2];
\draw[black, fill=mygreylight] (1.07021, -0.24436) circle [radius=0.125*0.2];
\draw[black, fill=mygreylight] (0.958894, -0.0939147) circle [radius=0.125*0.2];
\draw[black, fill=mygreylight] (0.90222, 0.0672621) circle [radius=0.125*0.2];
\draw[black, fill=mygreylight] (0.635158, 0.269538) circle [radius=0.125*0.2];
\draw[black, fill=mygreylight] (0.531446, 0.280634) circle [radius=0.125*0.2];
\draw[black, fill=mygreylight] (0.60217, 0.248038) circle [radius=0.125*0.2];
\draw[black, fill=mygreylight] (0.491523, 0.267716) circle [radius=0.125*0.2];
\draw[black, fill=mygreylight] (0.337463, 0.434121) circle [radius=0.125*0.2];
\draw[black, fill=mygreylight] (0.262031, 0.595759) circle [radius=0.125*0.2];
\draw[black, fill=mygreylight] (0.210778, 0.553144) circle [radius=0.125*0.2];
\draw[black, fill=mygreylight] (0.125937, 0.580563) circle [radius=0.125*0.2];
\draw[black, fill=mygreylight] (0.0377875, 0.708622) circle [radius=0.125*0.2];
\draw[black, fill=mygreylight] (-0.137369, 0.990575) circle [radius=0.125*0.2];
\draw[black, fill=mygreylight] (-0.0658768, 0.872351) circle [radius=0.125*0.2];
\draw[black, fill=mygreylight] (-0.193578, 0.84988) circle [radius=0.125*0.2];
\draw[black, fill=mygreylight] (1.5117, -0.543633) circle [radius=0.125*0.2];
\draw[black, fill=mygreylight] (1.24287, -0.211575) circle [radius=0.125*0.2];
\draw[black, fill=mygreylight] (1.19837, -0.161122) circle [radius=0.125*0.2];
\draw[black, fill=mygreylight] (1.01817, -0.087827) circle [radius=0.125*0.2];
\draw[black, fill=mygreylight] (0.777862, 0.226056) circle [radius=0.125*0.2];
\draw[black, fill=mygreylight] (0.804225, 0.147044) circle [radius=0.125*0.2];
\draw[black, fill=mygreylight] (0.621698, 0.428757) circle [radius=0.125*0.2];
\draw[black, fill=mygreylight] (0.55767, 0.354138) circle [radius=0.125*0.2];
\draw[black, fill=mygreylight] (0.483028, 0.542507) circle [radius=0.125*0.2];
\draw[black, fill=mygreylight] (0.384917, 0.595927) circle [radius=0.125*0.2];
\draw[black, fill=mygreylight] (0.407934, 0.614786) circle [radius=0.125*0.2];
\draw[black, fill=mygreylight] (0.222381, 0.660524) circle [radius=0.125*0.2];
\draw[black, fill=mygreylight] (0.118651, 0.673062) circle [radius=0.125*0.2];
\draw[black, fill=mygreylight] (0.00279014, 0.98281) circle [radius=0.125*0.2];
\draw[black, fill=mygreylight] (-0.156847, 1.08041) circle [radius=0.125*0.2];
\draw[black, fill=mygreylight] (-0.0952384, 0.984572) circle [radius=0.125*0.2];
\draw[black, fill=mygreylight] (1.50461, -0.463839) circle [radius=0.125*0.2];
\draw[black, fill=mygreylight] (1.73698, -0.569775) circle [radius=0.125*0.2];
\draw[black, fill=mygreylight] (1.42791, -0.255073) circle [radius=0.125*0.2];
\draw[black, fill=mygreylight] (1.11912, -0.0157329) circle [radius=0.125*0.2];
\draw[black, fill=mygreylight] (0.945205, 0.223911) circle [radius=0.125*0.2];
\draw[black, fill=mygreylight] (0.794893, 0.184432) circle [radius=0.125*0.2];
\draw[black, fill=mygreylight] (0.949975, 0.233761) circle [radius=0.125*0.2];
\draw[black, fill=mygreylight] (0.730641, 0.31773) circle [radius=0.125*0.2];
\draw[black, fill=mygreylight] (0.584427, 0.510088) circle [radius=0.125*0.2];
\draw[black, fill=mygreylight] (0.518569, 0.468361) circle [radius=0.125*0.2];
\draw[black, fill=mygreylight] (0.367593, 0.616443) circle [radius=0.125*0.2];
\draw[black, fill=mygreylight] (0.430644, 0.76937) circle [radius=0.125*0.2];
\draw[black, fill=mygreylight] (0.162015, 0.912671) circle [radius=0.125*0.2];
\draw[black, fill=mygreylight] (0.153808, 0.800495) circle [radius=0.125*0.2];
\draw[black, fill=mygreylight] (0.0887575, 1.1156) circle [radius=0.125*0.2];
\draw[black, fill=mygreylight] (-0.146548, 1.07395) circle [radius=0.125*0.2];
\draw[black, fill=mygreylight] (-0.149598, 1.18656) circle [radius=0.125*0.2];
\draw[black, fill=mygreylight] (-0.304762, 1.32004) circle [radius=0.125*0.2];
\draw[black, fill=mygreylight] (-0.365547, 1.59795) circle [radius=0.125*0.2];
\draw[black, fill=mygreylight] (2.66028, -1.50058) circle [radius=0.125*0.2];
\draw[black, fill=mygreylight] (1.22839, -0.101645) circle [radius=0.125*0.2];
\draw[black, fill=mygreylight] (1.16578, -0.0124005) circle [radius=0.125*0.2];
\draw[black, fill=mygreylight] (1.06987, 0.230152) circle [radius=0.125*0.2];
\draw[black, fill=mygreylight] (0.809648, 0.442526) circle [radius=0.125*0.2];
\draw[black, fill=mygreylight] (0.622528, 0.510586) circle [radius=0.125*0.2];
\draw[black, fill=mygreylight] (0.598886, 0.579151) circle [radius=0.125*0.2];
\draw[black, fill=mygreylight] (0.544, 0.588243) circle [radius=0.125*0.2];
\draw[black, fill=mygreylight] (0.53484, 0.652374) circle [radius=0.125*0.2];
\draw[black, fill=mygreylight] (0.456341, 0.862316) circle [radius=0.125*0.2];
\draw[black, fill=mygreylight] (0.25951, 1.06004) circle [radius=0.125*0.2];
\draw[black, fill=mygreylight] (0.141472, 1.1651) circle [radius=0.125*0.2];
\draw[black, fill=mygreylight] (0.0965729, 1.10694) circle [radius=0.125*0.2];
\draw[black, fill=mygreylight] (-0.115888, 1.29915) circle [radius=0.125*0.2];
\draw[black, fill=mygreylight] (-0.359633, 1.66913) circle [radius=0.125*0.2];
\draw[black, fill=mygreylight] (-0.384355, 1.68888) circle [radius=0.125*0.2];
\draw[black, fill=mygreylight] (-0.508018, 2.14324) circle [radius=0.125*0.2];
\draw[black, fill=mygreylight] (1.39521, -0.0818683) circle [radius=0.125*0.2];
\draw[black, fill=mygreylight] (1.44512, -0.0352168) circle [radius=0.125*0.2];
\draw[black, fill=mygreylight] (1.08602, 0.271529) circle [radius=0.125*0.2];
\draw[black, fill=mygreylight] (1.064, 0.312547) circle [radius=0.125*0.2];
\draw[black, fill=mygreylight] (0.833988, 0.582846) circle [radius=0.125*0.2];
\draw[black, fill=mygreylight] (0.833488, 0.414929) circle [radius=0.125*0.2];
\draw[black, fill=mygreylight] (0.634613, 0.591548) circle [radius=0.125*0.2];
\draw[black, fill=mygreylight] (0.59428, 0.716199) circle [radius=0.125*0.2];
\draw[black, fill=mygreylight] (0.63817, 0.811501) circle [radius=0.125*0.2];
\draw[black, fill=mygreylight] (0.450338, 0.862575) circle [radius=0.125*0.2];
\draw[black, fill=mygreylight] (0.405285, 1.08375) circle [radius=0.125*0.2];
\draw[black, fill=mygreylight] (0.258131, 1.12271) circle [radius=0.125*0.2];
\draw[black, fill=mygreylight] (0.154316, 1.08677) circle [radius=0.125*0.2];
\draw[black, fill=mygreylight] (-0.216573, 1.70811) circle [radius=0.125*0.2];
\draw[black, fill=mygreylight] (-0.198271, 1.77895) circle [radius=0.125*0.2];
\draw[black, fill=mygreylight] (-0.353702, 1.79722) circle [radius=0.125*0.2];
\draw[black, fill=mygreylight] (1.77775, -0.42834) circle [radius=0.125*0.2];
\draw[black, fill=mygreylight] (2.26136, -1.9869) circle [radius=0.125*0.2];
\draw[black, fill=mygreylight] (2.12801, -2.05305) circle [radius=0.125*0.2];
\draw[black, fill=mygreylight] (1.25128, 0.188517) circle [radius=0.125*0.2];
\draw[black, fill=mygreylight] (1.17953, 0.251889) circle [radius=0.125*0.2];
\draw[black, fill=mygreylight] (0.89172, 0.517588) circle [radius=0.125*0.2];
\draw[black, fill=mygreylight] (0.904357, 0.573754) circle [radius=0.125*0.2];
\draw[black, fill=mygreylight] (0.730123, 0.626772) circle [radius=0.125*0.2];
\draw[black, fill=mygreylight] (0.588188, 0.890863) circle [radius=0.125*0.2];
\draw[black, fill=mygreylight] (0.55663, 0.957573) circle [radius=0.125*0.2];
\draw[black, fill=mygreylight] (0.51405, 0.98785) circle [radius=0.125*0.2];
\draw[black, fill=mygreylight] (0.36083, 1.18816) circle [radius=0.125*0.2];
\draw[black, fill=mygreylight] (0.236824, 1.14339) circle [radius=0.125*0.2];
\draw[black, fill=mygreylight] (0.688785, -1.75172) circle [radius=0.125*0.2];
\draw[black, fill=mygreylight] (0.572976, -1.51133) circle [radius=0.125*0.2];
\draw[black, fill=mygreylight] (0.447526, -1.46017) circle [radius=0.125*0.2];

	\end{axis}
    \begin{pgfonlayer}{background}
    \draw[fill=white, drop shadow]
        (myaxis.north west) 
        rectangle (myaxis.south east);
    \end{pgfonlayer}
\end{tikzpicture}